\documentclass[11pt]{amsart}
\usepackage[T1]{fontenc}
\usepackage{lmodern}
\usepackage{microtype}
\usepackage[margin=1.12in]{geometry}
\usepackage{amsmath,amssymb,amsthm,mathtools}
\usepackage{booktabs,array}
\usepackage{enumitem}
\usepackage{needspace}
\usepackage{etoolbox}
\usepackage{tcolorbox}
\tcbuselibrary{skins,breakable}
\usepackage[hidelinks,pdfusetitle]{hyperref}

\makeatletter
\renewcommand\subsection{\@startsection{subsection}{2}{0pt}%
  {1.6ex plus .3ex}{.65ex}{\normalfont\bfseries}}
\makeatother
\AtBeginDocument{%
  \setlength{\parindent}{0pt}%
  \setlength{\parskip}{4.5pt plus 1pt minus 1pt}%
  \setlength{\jot}{5pt}%
  \raggedbottom
}

\setlist[enumerate]{itemsep=3.5pt, topsep=4pt}
\setlist[itemize]{itemsep=3.5pt, topsep=4pt}
\allowdisplaybreaks[1]

\newtcolorbox{mainbox}{
  enhanced,
  breakable,
  frame hidden,
  borderline west={1.6pt}{0pt}{black!60},
  colback=black!2,
  left=9pt,
  right=7pt,
  top=5pt,
  bottom=5pt,
  before skip=8pt,
  after skip=8pt
}

\newtheorem{theorem}{Theorem}[section]
\newtheorem{proposition}[theorem]{Proposition}
\newtheorem{lemma}[theorem]{Lemma}
\newtheorem{corollary}[theorem]{Corollary}
\newtheorem{conjecture}[theorem]{Conjecture}
\newtheorem*{theoremA}{Theorem A}
\newtheorem*{theoremB}{Theorem B}
\newtheorem*{theoremC}{Theorem C}
\newtheorem*{theoremD}{Theorem D}
\newtheorem*{theoremE}{Theorem E}
\newtheorem*{conjectureF}{Conjecture F}
\theoremstyle{definition}

\theoremstyle{remark}
\newtheorem{remark}[theorem]{Remark}

\numberwithin{equation}{section}

\newcommand{\Q}{\mathbb{Q}}
\newcommand{\Z}{\mathbb{Z}}
\newcommand{\C}{\mathbb{C}}
\newcommand{\R}{\mathbb{R}}
\newcommand{\HH}{\mathbb{H}}
\newcommand{\cO}{\mathcal{O}}
\newcommand{\cR}{\mathcal{R}}
\newcommand{\cP}{\mathcal{P}}
\newcommand{\cD}{\mathcal{D}}
\newcommand{\cU}{\mathcal{U}}
\newcommand{\Cinf}{\mathcal{C}_\infty}
\newcommand{\Cl}{\mathrm{Cl}}
\newcommand{\Nm}{\mathrm{N}}
\newcommand{\res}{\operatorname{res}}
\newcommand{\crr}{\operatorname{cr}}
\newcommand{\Gal}{\operatorname{Gal}}
\newcommand{\Ind}{\operatorname{Ind}}
\newcommand{\covol}{\operatorname{covol}}
\newcommand{\Disc}{\operatorname{Disc}}
\newcommand{\Rea}{\operatorname{Re}}
\newcommand{\Ima}{\operatorname{Im}}

\begin{document}

\title[Mahler measures and a cubic Dedekind zeta value]{Mahler measures and a cubic Dedekind zeta value in three dilogarithms}
\hypersetup{pdftitle={Mahler measures and a cubic Dedekind zeta value in three dilogarithms}}

\author{Vedran Men\dj u\v si\'c}
\address{Independent researcher, Croatia}
\date{19 November 2026}


\begin{abstract}
This paper proves an explicit evaluation formula for the three-variable Mahler measure
$n_4(s)=4\,m(x^4+y^4+z^4+1+s^{1/4}xyz)$ at the Heegner points
$\tau_d=(1+\sqrt{-d})/2$ for every fundamental discriminant $-d$ with $d>3$ and
$d\equiv3\pmod8$. Because $2$ is inert in $K=\Q(\sqrt{-d})$, the order of conductor $2$
has ring class number $3h(K)$, introducing ray class characters modulo $2$ that lie
beyond the reach of genus theory. The formula expresses $n_4(s_4(\tau_d))$ in terms
of complete traces of weight-three CM newforms, class-group twists of $\zeta_K(2)$,
and weight-one theta series attached to ray class characters modulo $2$. Every
individual root number occurring in these traces equals $+1$, and no splitting of
the ray class group or order-three assumption on its characters is required.

When $h(K)=1$, the ray-class contribution simplifies to $\zeta_F'(-1)$ for the cubic
field $F$ generated by the parameter $s_4(\tau_d)$. Combined with a separate primitive
normalization at the non-fundamental discriminant $-27$, this proves the six numerical
identities proposed by Zheng for $d\in\{11,19,27,43,67,163\}$.

In the first class-number-two case, $d=35$, the scalar part involves the cubic
field $C=\Q(a)$ ($a^3+2a-2=0$) and its quadratic extension $E=C(\sqrt5)$. An
exact hyperbolic volume-cycle certificate proves that $\zeta_C'(-1)=D(\Xi_C)/(6\pi)$
for an explicit three-term element $\Xi_C$ of the Bloch group. The analogous sextic
identity for $E$ remains conjectural and is reduced to a single explicit equality
in the Bloch group of the ray class field. Throughout, floating-point comparisons are
kept strictly separate from the exact algebraic and interval certificates used in the proofs.
\end{abstract}

\maketitle
\tableofcontents

\section{Introduction}\label{sec:intro}

\let\CMcompactsubsection\subsection
\pretocmd{\section}{\Needspace{170pt}}{}{}
\pretocmd{\subsection}{\Needspace{145pt}}{}{}
\BeforeBeginEnvironment{theorem}{\Needspace{155pt}}
\BeforeBeginEnvironment{proposition}{\Needspace{155pt}}
\BeforeBeginEnvironment{lemma}{\Needspace{155pt}}
\BeforeBeginEnvironment{corollary}{\Needspace{155pt}}
\BeforeBeginEnvironment{conjecture}{\Needspace{145pt}}
\BeforeBeginEnvironment{remark}{\Needspace{120pt}}
\BeforeBeginEnvironment{proof}{\Needspace{120pt}}
\clubpenalty=10000
\widowpenalty=10000

\subsection{The question}
For a Laurent polynomial $P$ in three variables, its logarithmic Mahler measure is
\[
m(P)=\int_{[0,1]^3}\log\bigl|P(e^{2\pi iu},e^{2\pi iv},e^{2\pi iw})\bigr|\,du\,dv\,dw.
\]
Following Samart, consider the one-parameter family
\begin{equation}\label{eq:n4}
n_4(s)=4\,m\bigl(x^4+y^4+z^4+1+s^{1/4}xyz\bigr).
\end{equation}
The choice of fourth root does not affect the integral, since replacing $x$ by $ix$
preserves Haar measure on the unit torus. Building on the modular parametrization of
Rogers and Samart \cite{Rogers09,Samart13,Samart15}, which expresses $n_4$ through
Eisenstein--Kronecker series, many evaluations at CM points have been written as rational
linear combinations of $L$-values of CM newforms and Dirichlet $L$-functions
\cite{GPQ21,HeYe22,Fei21}.

In those cases, genus theory converts the relevant Epstein zeta values at $s=2$ into
products of Dirichlet $L$-values \cite[Theorem~3.1]{HeYe22}; a general modular treatment
for elementary $2$-class groups is given in \cite[Corollary~5.4]{Fei21}. As Zheng
observed \cite[Section~9]{Zheng26}, this mechanism changes at the Heegner points
\[
\tau_d=\frac{1+\sqrt{-d}}2,\qquad d\equiv3\pmod8.
\]
For fundamental $-d$ with $d>3$, $2$ is inert in $K=\Q(\sqrt{-d})$, and the order
of conductor $2$ has ring class
number $3h(K)$. The new characters have nontrivial restriction to the residual
subgroup of order three; when $h(K)=1$, they are cubic. Their $L$-values lie outside
Samart's list of constants. Guided by $60$-digit numerical identities for
$d\in\{11,19,27,43,67,163\}$, Zheng formulated a general conjecture expressing such
Mahler measures in terms of weight-three CM forms and ray class $L$-values, while
leaving open the precise rule governing the levels and conductors
\cite[Conjectures~9.2 and~9.6]{Zheng26}.

This paper determines that rule and proves the formula for every fundamental
$d>3$ with $d\equiv3\pmod8$. The ray-class term has a concrete interpretation: when $h(K)=1$, it is the derivative at $s=-1$ of
the Dedekind zeta function of the cubic field generated by the parameter $s_d$ itself.
In the first class-number-two case, $d=35$, the ray-class contribution splits into a
cubic and a sextic piece, leading naturally to explicit Bloch--Wigner dilogarithms.

\subsection{Notation}
Let $\eta$ be the Dedekind eta function, $q=e^{2\pi i\tau}$, and set
\begin{equation}\label{eq:t-s4}
t(\tau)=\Bigl(\frac{\eta(\tau)}{\eta(2\tau)}\Bigr)^{24},\qquad
s_4(\tau)=\frac{(t(\tau)+64)^2}{t(\tau)},\qquad s_d=s_4(\tau_d).
\end{equation}
This is Samart's parameter for \eqref{eq:n4}; see Lemma~\ref{lem:s4}.
For $K=\Q(\sqrt{-d})$ with $d>3$, $d\equiv3\pmod8$ and $-d$ fundamental, let
$\cO=\Z+\tau_d\Z$ be its ring of integers and $h=h(K)$. Let $\Cl_{(2)}(K)$ denote the
ray class group modulo $(2)$, which has order $3h$ (Section~\ref{ssec:ray}). Write
\begin{itemize}
\item $\widehat{\Cl}$ for the $h$ characters of the class group, viewed as ray characters;
\item $\cR$ for the remaining $2h$ ray characters, all of exact conductor $(2)$;
\item $\cP\subset\cR$ for a set of representatives of the inverse pairs $\{\chi,\chi^{-1}\}$.
\end{itemize}
Fix an unramified Hecke character $\psi$ of $K$ of infinity type $2$, so that
$\psi((\alpha))=\alpha^2$. For a Hecke character $\Psi$, write $g_\Psi$ for its theta
series, and for a ray character $\chi$, set
$h_\chi=\sum_{(\mathfrak a,2)=1}\chi(\mathfrak a)q^{\Nm\mathfrak a}$. All $L$-functions
of $q$-series follow the convention $L(\sum a_nq^n,u)=\sum a_nn^{-u}$.

\subsection{Results}

\begin{mainbox}
\begin{theoremA}[Theorem~\ref{thm:GT}]
Let $d>3$ with $d\equiv3\pmod8$ and $-d$ a fundamental discriminant, and let $h=h(K)$.
Then
\[
\begin{aligned}
n_4(s_d)={}&\frac{40}{dh}\sum_{\eta\in\widehat{\Cl}}L'(g_{\psi\eta},0)
+\frac{40}{3dh}\sum_{\chi\in\cR}L'(g_{\psi\chi},0)\\[4pt]
&+\frac{5\sqrt d}{h\pi^3}\sum_{\eta\in\widehat{\Cl}}L_K(\eta,2)
-\frac{80}{3dh}\sum_{\rho\in\cP}L'(h_\rho,-1).
\end{aligned}
\]
The three traces satisfy functional equations with sign $+1$, and every individual
root number occurring in them equals $+1$ (Proposition~\ref{prop:signs}).
\end{theoremA}
\end{mainbox}

This formula requires neither a splitting $\Cl_{(2)}(K)\cong C_3\times\Cl(K)$ nor the
assumption that every character in $\cR$ has order three. In particular, class numbers
divisible by $3$ are covered without change, provided characters are paired with their
inverses rather than their squares.

When $h(K)=1$, the set $\cR=\{\chi_d,\chi_d^{-1}\}$ consists of the two cubic ray
characters. Put
\begin{equation}\label{eq:consts}
M_d=L'(g_{\psi},0),\qquad U_d=L'(g_{\psi\chi_d},0)+L'(g_{\psi\chi_d^{-1}},0),\qquad
\Delta_a=L'(\chi_{-a},-1),
\end{equation}
where $\chi_{-a}$ is the Kronecker character of $\Q(\sqrt{-a})$.

\begin{mainbox}
\begin{theoremB}[Theorem~\ref{thm:six}]
For $d\in\{11,19,43,67,163\}$ let $F_d=\Q(s_d)$. Then $F_d$ is a cubic field of
discriminant $-4d$, the minimal polynomial of $s_d$ is the cubic \eqref{eq:Qj}, and
\[
n_4(s_d)=\frac{10}{9d}\Bigl(36M_d+12U_d+3\Delta_d+288\,\zeta_{F_d}'(-1)\Bigr).
\]
\end{theoremB}
\end{mainbox}

The constant $\cD_d=36\zeta_{F_d}'(-1)$ accounts for the new ray-class term in the
numerical identities of \cite[Section~9]{Zheng26}. Arithmetically, it arises from a
two-dimensional irreducible Artin representation rather than a genus character. No claim is made regarding the linear independence of its special value from Dirichlet $L$-values.

\begin{mainbox}
\begin{theoremC}[Theorem~\ref{thm:d27}]
For $d=27$, $F_{27}=\Q(s_{27})\cong\Q(\sqrt[3]2)$ and
\[
n_4(s_{27})=\frac{10}{81}\Bigl(12M_{27}+4M_{108}+124M_{12}+25\Delta_3+96\,\zeta'_{\Q(\sqrt[3]2)}(-1)\Bigr),
\]
where $M_N=L'(g_N,0)$ for the primitive weight-three newforms $g_N$ of level $N$
constructed in Section~\ref{sec:d27}. The lattice constant in this case equals
$M_{108}+31M_{12}$ (Proposition~\ref{prop:F-lat}).
\end{theoremC}
\end{mainbox}

For $d=35$, the class number is $h=2$. Let $C=\Q(a)$ with $a^3+2a-2=0$ (discriminant
$-140$, one complex place), let $E=C(\sqrt5)$ (discriminant $98\,000$, signature $(2,2)$),
and let $\chi$ be the quadratic character of $E/C$. Denote by $\rho$ the two-dimensional
Artin representation attached to $C$ and by $\chi_5$ the Dirichlet character of $\Q(\sqrt5)$.

\begin{mainbox}
\begin{theoremD}[Theorem~\ref{thm:T35}]
Let $M_\pm$ be the derivatives at $0$ of the two weight-three CM newforms of level $35$,
and let $\cU_{140}$ be the sum of $L'(\cdot,0)$ over the four newforms of level $140$. Then
\[
n_4(s_{35})=\frac47(M_++M_-)+\frac4{21}\cU_{140}+\frac1{21}\Delta_{35}+\frac8{35}\Delta_7
+\frac{32}{7}\zeta_C'(-1)+\frac{20}{21}L'_{E/C}(-1,\chi).
\]
Moreover, $L(s,\chi)=L(s,\chi_5)\,L(s,\rho\otimes\chi_5)$ and $L(-1,\chi_5)=-\tfrac25$.
\end{theoremD}
\end{mainbox}

The expressions below use the Bloch--Wigner dilogarithm $D$. The existence of a
dilogarithmic formula for $\zeta_C'(-1)$ is known \cite{Borel77,Zagier86,Suslin91};
the contribution here is the explicit three-term class and its coefficient.
The sextic identity remains conjectural.
Set
\[
\Xi_C=[-a^5]-5[a]-15[-a].
\]
With $r=\sqrt5$, $\sigma\colon r\mapsto-r$ fixing $a$, and $\eta_\sigma(z)=[z]-[\sigma z]$, put
\[
\Xi_E=2\eta_\sigma(\beta^{-1})-\tfrac23\eta_\sigma\Bigl(-\frac1{2u}\Bigr)
-3\eta_\sigma\Bigl(-\frac1{\beta x}\Bigr)+2\eta_\sigma(x^{-2})+\eta_\sigma\Bigl(-\frac a{2\beta u}\Bigr),
\]
where $\beta=r/(a+1)$, $x=(\beta-1)/a$, and $u=(1+r)/2$. Throughout, values of $D$ are
evaluated at the complex embedding with $\Ima a>0$ and $r>0$.

\begin{mainbox}
\begin{theoremE}[Theorem~\ref{thm:R1}]
$\zeta_C'(-1)=\frac{D(\Xi_C)}{6\pi}$.
\end{theoremE}
\end{mainbox}

\begin{mainbox}
\begin{conjectureF}[Conjecture~\ref{conj:R2}]
$L'_{E/C}(-1,\chi)=\frac{6}{5\pi}D(\Xi_E)$; equivalently,
\[
L'(-1,\rho\otimes\chi_5)=-\frac3\pi D(\Xi_E)
\qquad\text{or}\qquad
L(2,\rho\otimes\chi_5)=\frac{24\pi^2}{140^{3/2}}D(\Xi_E).
\]
\end{conjectureF}
\end{mainbox}

At $130$-digit working precision, the observed residual for Conjecture~F is
approximately $4.51\cdot10^{-130}$; this is numerical evidence rather than a certified
error bound. By Sato's theorem \cite{Sato26}, Conjecture~F is equivalent to the identity
$A_6=-3\Xi_E$ inside a one-dimensional $\Q$-vector space (Corollary~\ref{cor:A6}), while
Proposition~\ref{prop:noKlein} explains why the hyperbolic volume argument used for $C$
cannot be applied directly to $E$.

Substituting Theorem~E into Theorem~D expresses the cubic contribution in closed
dilogarithmic form:
\[
n_4(s_{35})=\frac47(M_++M_-)+\frac4{21}\cU_{140}+\frac1{21}\Delta_{35}+\frac8{35}\Delta_7
+\frac{16}{21\pi}D(\Xi_C)+\frac{20}{21}L'_{E/C}(-1,\chi).
\]
Conjecture~F would turn the final term into $\frac{8}{7\pi}D(\Xi_E)$.

\subsection{Relation to earlier work}\label{ssec:literature}
On the analytic side, the starting point is Samart's Eisenstein--Kronecker lattice
formula \eqref{eq:lattice} \cite[Proposition~2.1]{Samart13}, which He and Ye state for
all $\tau$ with $\Ima\tau\ge1/\sqrt2$ alongside its CM theta--Epstein specialization
\cite[Propositions~5.2--5.3, (5.2)]{HeYe22}. At $\tau_d$, that specialization corresponds
to the quadratic form $a=1$, $b=-1$, $c=(1+d)/4$, and every point considered here
satisfies $\Ima\tau_d\ge\sqrt{11}/2>1/\sqrt2$.

He and Ye's Table~1 also lists these six rational CM $j$-invariants for the family $n_2$.
If $S$ is their parameter satisfying $j=(S-16)^3/S$, then setting $S=-4096/t$ yields
\[
s=128-S-\frac{4096}{S},\qquad
S=\frac{s^2+(513-j)s+127j+107136}{j+3375}.
\]
Because $j\ne-3375$ for all six values here, the parameter fields $\Q(s_d)$ and
$\Q(S_d)$ coincide; note, however, that this is purely an algebraic relation between
parameters and not an identity between $n_2$ and $n_4$.

The arithmetic contributions are the decomposition into complete ray class character
traces modulo $(2)$ and their root numbers (Theorem~A), the identification of the
ray-class term with cubic parameter fields, the primitive normalization at $-27$, and
the analysis of $d=35$ together with its dilogarithms. These results extend the
genus-character calculations of \cite{HeYe22}, the elementary-$2$-class-group
framework of \cite[Corollary~5.4]{Fei21}, and the coset decompositions of
\cite[Section~9]{Zheng26}. Theorems~B and~C prove the six identities
proposed in \cite[Section~9.1]{Zheng26}, using standard analytic continuation and
imprimitive Euler-factor techniques \cite{GPQ21,Fei21}.

\subsection{The domain hypothesis at \texorpdfstring{$\tau_d$}{tau d}}\label{ssec:domain-lit}
For the formula used here, Samart's Proposition~2.1 assumes a positive real nome $q\in(0,1)$ as well as
$\Ima\tau\ge1/\sqrt2$ \cite[p.~5]{Samart13}. At $\tau_d$, the nome is
$q=-e^{-\pi\sqrt d}<0$, so an analytic continuation is required before the lattice
formula can be applied.

Zheng's Lemma~2.13 continues the identity on the connected component $\Cinf$ of
$\{|s_4|>256\}$ containing a punctured neighbourhood of the cusp \cite{Zheng26}. A
pointwise interval check verifies the inequality $|s_4(\tau_d)|>256$, but not
membership in $\Cinf$; his Remark~9.4 reports a $40$-digit numerical scan along the
vertical ray from $\tau_d$ without an interval certification of the entire ray.

Lemma~\ref{lem:domain} resolves this uniformly by showing that $|s_4(\tau)|>847$ on the
entire connected half-plane $\Ima\tau\ge3/2$, placing every $\tau_d$ ($d\ge11$) inside
$\Cinf$. Proposition~\ref{prop:continuation} recalls Zheng's continuation lemma
with a short proof from Rogers's local identity.

Finally, the minimal polynomials of $s_d$ were left open in \cite[Section~9.1]{Zheng26}
because a $110$-digit integer-relation search ruled out cubics with coefficients below
$10^8$. Lemma~\ref{lem:s4} and Proposition~\ref{prop:minpoly} supply the exact minimal
cubics: their constant terms $(j_d-1728)^2$ exceed $10^8$ in every case, starting at
$1\,189\,974\,016$ for $d=11$.

\Needspace{80pt}
\CMcompactsubsection{Organization}
Sections~\ref{sec:parameters}--\ref{sec:general} develop the modular parameters, the
analytic continuation to lattice sums, and the general trace formula.
Sections~\ref{sec:classone}--\ref{sec:d35} carry out the explicit evaluations for
$h(K)=1$, $d=27$, and $d=35$. Section~\ref{sec:dilog} treats the dilogarithms at
$d=35$, and Section~\ref{sec:computations} summarizes the exact, interval, and
numerical computations supporting the paper.

\section{The parameters and the cusp component}\label{sec:parameters}

Put $W(\tau)=\eta(\tau)\eta(4\tau)^2/\eta(2\tau)^3$. Samart's parameter for
\eqref{eq:n4} is $t^{-1}(16W^4+W^{-4})^4$ \cite[Proposition~2.1]{Samart13}, which is also
used in \cite[Proposition~5.2]{HeYe22}; it is treated here as a holomorphic function of $\tau$.

\begin{lemma}\label{lem:s4}
As meromorphic functions on $\HH$,
\begin{equation}\label{eq:s4-identity}
t^{-1}\bigl(16W^4+W^{-4}\bigr)^4=\frac{(t+64)^2}{t}=s_4 .
\end{equation}
If $j=j(\tau)$ and $s=s_4(\tau)$, then $Q_j(s)=0$, where
\begin{equation}\label{eq:Qj}
Q_j(s)=s^3+(432-j)s^2+(207j+62208)s+(j-1728)^2,
\end{equation}
and $\operatorname{disc}_sQ_j=4j^2(j-1728)(j+3375)^2$.
\end{lemma}

\begin{proof}
Write $A=\theta_2(2\tau)$, $B=\theta_3(2\tau)$, and $C=\theta_4(2\tau)$. The theta--eta identities
and Jacobi's quartic identity give
\[
A=\frac{2\eta(4\tau)^2}{\eta(2\tau)},\quad C=\frac{\eta(\tau)^2}{\eta(2\tau)},\quad
ABC=2\eta(2\tau)^3,\quad B^4=A^4+C^4 .
\]
Put $k=A^2/B^2$. Then $W^2=CA^2/(4\eta(2\tau)^3)=A/(2B)$, so $W^4=k/4$, and
\[
t=\frac{C^{12}}{\eta(2\tau)^{12}}=\frac{16(1-k^2)^2}{k^2}=\bigl(16W^4-W^{-4}\bigr)^2 .
\]
Since the squares of $16W^4\pm W^{-4}$ differ by $64$, it follows that $(16W^4+W^{-4})^2=t+64$,
which proves \eqref{eq:s4-identity} without extracting a square root.

Combining the level-two relation $j=(t+256)^3/t^2$ with the resultant identity
\[
\operatorname{Res}_t\bigl(t^2+(128-s)t+4096,\ (t+256)^3-jt^2\bigr)=2^{24}Q_j(s)
\]
shows that $Q_j(s)=0$. The discriminant follows by direct expansion.
\end{proof}

\begin{lemma}\label{lem:domain}
For every $\tau\in\HH$ with $\Ima\tau\ge3/2$,
\[
|t(\tau)|>\frac{325000}{333}>975,\qquad |s_4(\tau)|>847>256 .
\]
Consequently, every $\tau$ with $\Ima\tau\ge 3/2$ lies in the connected component
$\Cinf$ of $\{\tau\in\HH:|s_4(\tau)|>256\}$ containing a neighbourhood of the cusp
$i\infty$. In particular, $\tau_d\in\Cinf$ for every $d\ge11$, with lower bounds
independent of $d$.
\end{lemma}

\begin{proof}
From $t=q^{-1}\prod_{n\ge1}(1+q^n)^{-24}$ and $r=|q|\le e^{-3\pi}<10^{-3}$,
$\prod_{n\ge1}|1+q^n|\le\exp(r/(1-r))$, and hence
\[
|t|\ge r^{-1}\exp\Bigl(-\frac{24r}{1-r}\Bigr)>1000\Bigl(1-\frac{24}{999}\Bigr)=\frac{325000}{333}.
\]
Because $|t|>64$, it follows that $|s_4|\ge(|t|-64)^2/|t|=|t|-128+4096/|t|>847$.

The half-plane $\Ima\tau\ge3/2$ is connected, contains a neighbourhood of $i\infty$, and
is contained in $\{|s_4|>256\}$, so it lies entirely in $\Cinf$. Finally,
$\Ima\tau_d=\sqrt d/2\ge\sqrt{11}/2>3/2$.
\end{proof}

In particular, every vertical ray $\tau_d+iv$ ($v\ge0$) remains inside $\Cinf$ and tends
to the cusp, which is precisely what the analytic continuation in
Section~\ref{sec:analytic} requires (see Section~\ref{ssec:domain-lit}).

\clearpage
\begin{proposition}\label{prop:minpoly}
Let $d\in\{11,19,27,43,67,163\}$ and let $j_d$ be as in Table~\ref{tab:roots}. Then
$Q_{j_d}$ is the minimal polynomial of $s_d$ over $\Q$. It has a unique real root,
namely $s_d$, and this root lies in the interval $(-N_d,-N_d+1)$.
\end{proposition}

\begin{table}[ht]
\centering
\small
\setlength{\tabcolsep}{9pt}
\renewcommand{\arraystretch}{1.15}
\begin{tabular}{@{}rrrrc@{}}
\toprule
$d$ & $j_d$ & $N_d$ & $p_d$ & $(a,b,c)\bmod p_d$ \\
\midrule
$11$ & $-32^3$ & $33\,403$ & $3$ & $(2,0,1)$\\
$19$ & $-96^3$ & $96^3+640$ & $5$ & $(3,1,1)$\\
$27$ & $-12\,288\,000$ & $12\,288\,000+640$ & $7$ & $(2,4,2)$\\
$43$ & $-960^3$ & $960^3+640$ & $11$ & $(8,2,3)$\\
$67$ & $-5280^3$ & $5280^3+640$ & $17$ & $(4,14,13)$\\
$163$ & $-640\,320^3$ & $640\,320^3+640$ & $41$ & $(12,31,16)$\\
\bottomrule
\end{tabular}
\caption{Singular moduli, root intervals, and primes certifying irreducibility of
$Q_{j_d}\equiv x^3+ax^2+bx+c$.}
\label{tab:roots}
\end{table}

\begin{proof}
For $d\ne27$, the values $j_d=j(\tau_d)$ are the classical class-number-one singular
moduli \cite{YuiZagier97,Mendusic26}. For the order of discriminant $-27$, the CM
curve $y^2+y=x^3-30x+63$ \cite[Example~2.8]{Melistas24} has $c_4=1440$ and
$\Delta=-243$, giving $j=-12\,288\,000$; since $[1,1,7]$ is the unique reduced primitive
form of discriminant $-27$, this equals $j(\tau_{27})$. By Lemma~\ref{lem:s4}, $Q_{j_d}(s_d)=0$.

To verify irreducibility, let $Q\in\mathbb F_{p_d}[x]$ be the reduction listed in
Table~\ref{tab:roots}. The remainders of $x^{p_d}-x$ modulo $Q$ are, in row order,
\[
(r_2,r_1,r_0)=(1,2,2),\ (3,4,2),\ (6,5,2),\ (5,0,10),\ (12,3,1),\ (37,33,9).
\]
Their resultants with $Q$ are $1,2,1,5,12,18\ne0$, so $Q$ has no root in
$\mathbb F_{p_d}$ and is therefore irreducible. By Gauss's lemma, $Q_{j_d}$ is the
minimal polynomial of $s_d$.

Each $j_d$ is negative and distinct from $-3375$, so $\operatorname{disc}_s Q_{j_d}<0$ and
$Q_{j_d}$ has a unique real root. At $\tau_d$, $q=-e^{-\pi\sqrt d}$, so every
factor $1+q^n$ in the eta product is positive and $t(\tau_d)<0$. Lemma~\ref{lem:domain}
rules out $t=-64$, hence $s_d=(t+64)^2/t<0$ is real. Evaluating at integer endpoints
gives $Q_{j_d}(-N_d)<0<Q_{j_d}(-N_d+1)$ in every row.
\end{proof}

\begin{remark}\label{rem:fricke}
For $d\ne27$, if $P_d(X)=a_dX^3+b_dX^2+c_dX+4096$ is the cubic Fricke polynomial of
\cite[Section~12]{Mendusic26} in the variable $X=256t/(t+64)^2=256/s_4$, then
$Q_{j_d}(s)=\frac{s^3}{4096}P_d(256/s)$. At $d=11$, for instance,
\[
Q_{j_{11}}(s)=s^3+33\,200s^2-6\,720\,768s+1\,189\,974\,016,
\qquad -33\,403<s_{11}<-33\,402,
\]
and the corresponding Mahler measure $n_4(s_{11})\approx 10.41709\,66323\,51$ is
recorded in Table~\ref{tab:values}. At $d=163$, the same CM value governs the cubic
Ramanujan--Sato series of \cite{Mendusic16326} with parameter $X=256/s_{163}$, whose
generating function ${}_3F_2(\tfrac14,\tfrac12,\tfrac34;1,1;X)$ equals
$s\,d\mathcal M/ds$ at $s=s_{163}$ in the notation of Section~\ref{sec:analytic}.
\end{remark}

\begin{remark}\label{rem:general-d}
For general $d$ in the family, $s_d$ is still an algebraic CM value, though its
minimal polynomial plays no role in Theorem~A. For $d=35$, Lemma~\ref{lem:domain} applies
since $\Ima\tau_{35}=\sqrt{35}/2>3/2$; in the ancillary computations $s_{35}$ is a real
algebraic number of degree six.
\end{remark}

\section{From the Mahler measure to lattice sums}\label{sec:analytic}

\begin{lemma}[Rogers {\cite[Proposition~2.2, (2.9)]{Rogers09}}]\label{lem:local}
For $|s|>256$,
\[
n_4(s)=\Rea\mathcal M(s),\qquad
\mathcal M(s)=\log s-\sum_{n\ge1}\frac{(4n)!}{n(n!)^4}s^{-n},
\]
for any local branch of $\log s$. The differential
$d\mathcal M=\bigl(s^{-1}+\sum_{n\ge1}\frac{(4n)!}{(n!)^4}s^{-n-1}\bigr)ds$ is single valued
and holomorphic on $|s|>256$.
\end{lemma}

\begin{proof}
The constant-term argument follows the method of \cite{RV99}.
Choose $c$ with $c^4=s$ and put $Q=x^4+y^4+z^4+1$. On the unit torus, write
$Q+cxyz=cxyz(1+Q/(cxyz))$ with $|Q/(cxyz)|\le4/|c|<1$, so the logarithmic series of the
second factor converges uniformly. Its degree-$m$ term has a nonzero constant Fourier
coefficient only when $m=4n$, where that coefficient equals $(4n)!/((n!)^4c^{4n})$.
Multiplying $-1/(4n)$ by $4$ yields the series for $\mathcal M(s)$.

Normal convergence on $|s|>256$ follows from $(4n)!/(n!)^4\le256^n$, and replacing $c$
by $ic$ leaves the Mahler measure unchanged.
\end{proof}

Define a holomorphic function on $\HH$ by
\begin{equation}\label{eq:E4}
\mathcal E(\tau)=-2\pi i\tau-80\sum_{k\ge1}k^2\log(1-q^k)+160\sum_{k\ge1}k^2\log(1-q^{2k}),
\end{equation}
where each logarithm is given by its power series around $q=0$.

\begin{proposition}[Zheng {\cite[Lemma~2.13]{Zheng26}}]\label{prop:continuation}
For every $\tau\in\Cinf$, $n_4(s_4(\tau))=\Rea\mathcal E(\tau)$.
\end{proposition}

\begin{proof}
This uses the differential-comparison approach of \cite{ZhengInterior26,Zheng26}.
Rogers \cite[Theorem~2.3, (2.14)]{Rogers09} proves that for small $|q|$,
\[
n_4(s_4(q))=-\tfrac13G(q)+\tfrac23G(q^2),\qquad
G(q)=\Rea\Bigl(-\log q+240\sum_{k\ge1}k^2\log(1-q^k)\Bigr).
\]
For small positive $q$, $s_4(q)>256$, and both $\mathcal M(s_4(q))$
and the holomorphic lift of the right-hand side are real on that segment. Equality of
their real parts therefore implies equality of holomorphic functions, which extends
to a small disc by the identity theorem. Pulling back the right-hand side via
$q=e^{2\pi i\tau}$ gives precisely \eqref{eq:E4}.

Now the differentials $d(\mathcal M\circ s_4)$ and $d\mathcal E$ are single valued and
holomorphic on $\Cinf$. Since they agree on a nonempty open set, they agree throughout
the connected set $\Cinf$. Consequently, $n_4(s_4(\tau))-\Rea\mathcal E(\tau)$ is a
locally constant, single-valued real function on $\Cinf$; being zero on the initial
disc, it vanishes identically on $\Cinf$.
\end{proof}

The following standard series calculation also appears in
\cite[(5)--(6)]{Zheng26}. For each integer $j\ge1$, define
\[
U_j(\tau)=\sum_{m\ne0}\frac1m\sum_{n\in\Z}(jm\tau+n)^{-3}.
\]

\begin{lemma}\label{lem:Uj}
The series $U_j$ converges absolutely and locally uniformly on $\HH$, and
$U_j(\tau)=-8\pi^3i\sum_{k\ge1}k^2\log(1-q^{jk})$. Consequently,
\[
\mathcal E(\tau)=-i\Bigl(2\pi\tau+\frac{10}{\pi^3}\bigl(U_1(\tau)-2U_2(\tau)\bigr)\Bigr).
\]
\end{lemma}

\begin{proof}
For $\Ima\tau\ge y_0>0$ and $a=j|m|y_0$, the uniform bound
$\sum_{n\in\Z}((n+x)^2+a^2)^{-3/2}\le 2/a^3+2/a^2$ yields normal convergence. For
$\Ima z>0$, differentiating the partial-fraction expansion of $\pi^2\csc^2(\pi z)$ gives
$\sum_{n\in\Z}(z+n)^{-3}=4\pi^3i\sum_{k\ge1}k^2e^{2\pi ikz}$. Pairing $m<0$ with $m>0$
doubles the sum, giving $U_j=8\pi^3i\sum_{m,k\ge1}k^2q^{jmk}/m=-8\pi^3i\sum_{k\ge1}k^2\log(1-q^{jk})$.
\end{proof}

For a lattice $L\subset\C$, define the absolutely convergent sums
\begin{equation}\label{eq:BET}
A(L)=\sum_{\lambda\in L\setminus0}|\lambda|^{-4},\qquad
H(L)=\sum_{\lambda\in L\setminus0}\frac{\Rea(\lambda^2)}{|\lambda|^6},\qquad
T(L)=A(L)+2H(L).
\end{equation}

\begin{proposition}[Samart--He--Ye lattice formula]\label{prop:lattice}
For $\tau\in\Cinf$ and $L_j=\Z+j\tau\Z$,
\begin{equation}\label{eq:lattice}
n_4(s_4(\tau))=\frac{10\,\Ima\tau}{\pi^3}\bigl(-T(L_1)+4T(L_2)\bigr).
\end{equation}
\end{proposition}

\begin{proof}
Writing $\lambda=x+iv$, $\Ima(\lambda^{-3})=-v(3x^2-v^2)/|\lambda|^6$ and
$|\lambda|^{-4}+2\Rea(\lambda^2)|\lambda|^{-6}=(3x^2-v^2)|\lambda|^{-6}$. Setting
$\lambda=jm\tau+n$ and $y=\Ima\tau$ gives $v=jmy$, so
$\Ima U_j=-jy\bigl(T(L_j)-6\zeta(4)\bigr)$, where the row $m=0$ contributes
$\sum_{n\ne0}3n^{-4}=6\zeta(4)$. Taking real parts in Lemma~\ref{lem:Uj} yields
\[
\Rea\mathcal E=2\pi y+\frac{10y}{\pi^3}\Bigl(-T(L_1)+4T(L_2)-\frac{\pi^4}5\Bigr)
=\frac{10y}{\pi^3}\bigl(-T(L_1)+4T(L_2)\bigr),
\]
and Proposition~\ref{prop:continuation} completes the proof.
\end{proof}

\begin{remark}
Formula~\eqref{eq:lattice} matches the normalization of \cite[Proposition~2.1(iii)]{Samart13},
\cite[Proposition~5.2]{HeYe22}, and \cite[(7)]{Zheng26}. At $\tau=\tau_d$, the two lattices
are $\cO=\Z+\tau_d\Z$ and $\cO_2=\Z+2\cO$, so for every $d\ge11$,
\begin{equation}\label{eq:lattice-d}
n_4(s_d)=\frac{5\sqrt d}{\pi^3}\bigl(-T(\cO)+4T(\cO_2)\bigr).
\end{equation}
This is the case $a=1$, $b=-1$, $c=(1+d)/4$ of \cite[(5.2)]{HeYe22}. As an independent
numerical sanity check, direct lattice summation at $d=11$ and $d=35$ recovers the
Mahler-measure values of Table~\ref{tab:values} to within the truncation error.
\end{remark}

\section{The general formula}\label{sec:general}

Throughout this section assume $d>3$, $d\equiv3\pmod8$, and $-d$ fundamental, and
write $K=\Q(\sqrt{-d})$, $\cO=\Z+\tau_d\Z$, $h=h(K)$, and $\cO_2=\Z+2\cO$. Since $d>3$,
the unit group is $\cO^\times=\{\pm1\}$. Put $\Nm\alpha=\alpha\bar\alpha$.

\subsection{The ray class group modulo \texorpdfstring{$2$}{2}}\label{ssec:ray}
The minimal polynomial of $\tau_d$ is $u^2-u+(d+1)/4$, which reduces to $u^2+u+1$
modulo $2$. Hence $(2)$ is inert in $K$, $\cO/2\cO\cong\mathbb F_4$, and
$(\cO/2\cO)^\times=\langle\bar\tau_d\rangle$ is cyclic of order $3$. Since $\pm1\equiv1\pmod{2\cO}$,
the exact sequence
\[
1\to(\cO/2\cO)^\times\to\Cl_{(2)}(K)\to\Cl(K)\to1
\]
gives $|\Cl_{(2)}(K)|=3h$. Any ray character trivial on $(\cO/2\cO)^\times$ factors
through $\Cl(K)$, whereas every other ray character has conductor exactly $(2)$, as $(2)$
has no nontrivial proper ideal divisor. This accounts for the partition into $\widehat{\Cl}$,
$\cR$, and $\cP$ in the introduction, with $|\cR|=2h$ and $|\cP|=h$.

For any integral ideal $\mathfrak a$ coprime to $2$, the norm $\Nm\mathfrak a$ is odd.
Thus $\mathfrak a\bar{\mathfrak a}=(\Nm\mathfrak a)$ is generated by an element
congruent to $1$ modulo $2\cO$, so its class in $\Cl_{(2)}(K)$ is trivial. Because every
ray class contains such an ideal representative, complex conjugation acts by inversion on
$\Cl_{(2)}(K)$, giving $\chi(\bar{\mathfrak a})=\chi(\mathfrak a)^{-1}$ for every ray
character $\chi$. Since conjugation preserves norms,
\begin{equation}\label{eq:conj}
L_K(\chi,s)=L_K(\chi^{-1},s),\qquad h_\chi=h_{\chi^{-1}}.
\end{equation}
In general, the inverse of $\chi\in\cR$ need not equal $\chi^2$, so the pairing in
$\cP$ is always by inverses.

\subsection{Hecke characters and theta series}
Principal ideals form a finite-index subgroup of the group of fractional ideals, and
$(\alpha)\mapsto\alpha^2$ is well defined on it because $\cO^\times=\{\pm1\}$. By
divisibility of $\C^\times$, this map extends to all fractional ideals; fix one such
unramified Hecke character $\psi$ of infinity type $2$. Changing $\psi$ by a class-group
character merely permutes the sets $\{\psi\eta\}$ and $\{\psi\chi\}$ below.

For a Hecke character $\Psi$ of conductor $\mathfrak f$, write
$g_\Psi=\sum_{(\mathfrak a,\mathfrak f)=1}\Psi(\mathfrak a)q^{\Nm\mathfrak a}$. By
Hecke's theorem \cite[Theorem~3.7]{Fei21} (see also \cite{Kani14}), $g_{\psi\eta}$
($\eta\in\widehat{\Cl}$) is a weight-three cusp form of level $d$, while $g_{\psi\chi}$
($\chi\in\cR$) is a weight-three cusp form of level $4d$; both have nebentypus
$\chi_{-d}$ because their finite parts are trivial on rational integers prime to $2$.
Likewise, for $\chi\in\cR$, $h_\chi$ is a weight-one cusp form of level $4d$ and
nebentypus $\chi_{-d}$. Its cuspidality holds because $\chi\ne\chi^{-1}$: the
restriction of $\chi$ to $(\cO/2\cO)^\times$ has order three, so $\chi$ cannot be
quadratic regardless of the parity of $h$.

\subsection{Lattice sums and characters}
Define the character sums
\[
M_{\widehat{\Cl}}=\sum_{\eta\in\widehat{\Cl}}L(g_{\psi\eta},3),\quad
M_{\cR}=\sum_{\chi\in\cR}L(g_{\psi\chi},3),\quad
Z_{\widehat{\Cl}}=\sum_{\eta\in\widehat{\Cl}}L_K(\eta,2),\quad
Z_{\cP}=\sum_{\rho\in\cP}L_K(\rho,2).
\]

\begin{proposition}\label{prop:decomp}
With $T$ as in \eqref{eq:BET},
\[
-T(\cO)+4T(\cO_2)=\frac1h\Bigl(2M_{\widehat{\Cl}}+\frac{16}3M_{\cR}+Z_{\widehat{\Cl}}+\frac{16}3Z_{\cP}\Bigr).
\]
\end{proposition}

\begin{proof}
The calculation uses character orthogonality on $\cO$ and $1+2\cO$, removes the
Euler factor at $(2)$, and combines the two cosets of $\cO_2$.

\smallskip
\Needspace{90pt}
\noindent\textbf{The maximal-order sums.}\par
Let $H_0=H(\cO)$, $A_0=A(\cO)$, and let $H_1$, $A_1$ be the corresponding sums
restricted to $\alpha\equiv1\pmod{2\cO}$. Every nonzero principal ideal has two
generators, with $\alpha^2=\psi((\alpha))$, so character orthogonality on $\Cl(K)$ yields
\[
H_0=\frac2h M_{\widehat{\Cl}},\qquad A_0=\frac2hZ_{\widehat{\Cl}},
\]
where the harmonic sum is real because conjugation preserves $\cO$.

\smallskip
\Needspace{90pt}
\noindent\textbf{The ray-class sums.}\par
An ideal coprime to $2$ is trivial in $\Cl_{(2)}(K)$ if and only if it has a generator
$\equiv1\pmod{2\cO}$, in which case both generators satisfy this congruence since
$-1\equiv1\pmod{2\cO}$. Orthogonality on $\Cl_{(2)}(K)$ therefore gives
\[
H_1=\frac2{3h}\Bigl(\frac{15}{16}M_{\widehat{\Cl}}+M_{\cR}\Bigr),\qquad
A_1=\frac2{3h}\Bigl(\frac{15}{16}Z_{\widehat{\Cl}}+2Z_{\cP}\Bigr).
\]
Here the factor $15/16$ strips out the Euler factor at $(2)$ from the unramified
$L$-series: for infinity type $2$ at $s=3$ this factor is
$1-\psi((2))4^{-3}=1-4\cdot4^{-3}$, and for finite characters at $s=2$ it is
$1-4^{-2}$. The factor $2$ multiplying $Z_\cP$ comes from \eqref{eq:conj}.

\smallskip
\Needspace{90pt}
\noindent\textbf{Combining the cosets.}\par
Since $\cO_2=2\cO\sqcup(1+2\cO)$ and scaling by $2$ multiplies $A$ and $H$ by $1/16$,
$T(\cO_2)=T(\cO)/16+2H_1+A_1$. The stated expression follows by substitution in
\[
-T(\cO)+4T(\cO_2)=-\tfrac34(2H_0+A_0)+4(2H_1+A_1).\qedhere
\]
\end{proof}

\subsection{Complete traces and their functional equations}
For a lattice $L\subset\C$, set
\[
\Theta_L(z)=\frac12\sum_{\alpha\in L}\alpha^2q^{\Nm\alpha},\qquad
\vartheta_L(z)=\frac12\sum_{\alpha\in L}q^{\Nm\alpha}.
\]
On functions on the imaginary axis, define the normalized involution of weight $k$
and level $N$ by
\[
(W_{N,k}f)(iy)=N^{-k/2}y^{-k}f\bigl(i/(Ny)\bigr).
\]

\begin{lemma}\label{lem:poisson}
Let $L\subset\C$ be a lattice stable under complex conjugation, of covolume $A$, and
put $N=4A^2$. Then
\[
\Theta_L\Bigl(\frac i{Ny}\Bigr)=N^{3/2}y^3\,\Theta_L(iy),\qquad
\vartheta_L\Bigl(\frac i{Ny}\Bigr)=N^{1/2}y\,\vartheta_L(iy).
\]
\end{lemma}

\begin{proof}
Let $J$ denote rotation by $\pi/2$ and let $B$ be an oriented basis matrix for $L$ with
$\det B=A$. The identity $B^{-T}=A^{-1}JBJ^{-1}$ gives the Euclidean dual lattice
$L^\vee=(i/A)L$. Poisson summation applied to $e^{-2\pi y|x|^2}$ immediately gives
the second identity.

For the first identity, conjugation invariance implies
$\Theta_L(iy)=\frac12\sum_{\lambda\in L}\Rea(\lambda^2)e^{-2\pi y|\lambda|^2}$. With the
harmonic polynomial $P(x_1,x_2)=x_1^2-x_2^2$ and Fourier convention $e^{-2\pi ix\cdot\xi}$,
differentiating the Gaussian transform gives
\[
\widehat{Pe^{-2\pi y|x|^2}}(\xi)=-(2y)^{-3}P(\xi)e^{-\pi|\xi|^2/(2y)},
\]
as the constant terms of the two second derivatives cancel. Since
$P(i\lambda/A)=-A^{-2}P(\lambda)$, the two minus signs cancel, and Poisson summation
yields $\Theta_L(iy)=(8A^3y^3)^{-1}\,\Theta_L\bigl(i/(4A^2y)\bigr)$.
\end{proof}

The covolumes of $\cO$ and $\cO_2$ are $\sqrt d/2$ and $\sqrt d$. By
Lemma~\ref{lem:poisson},
\[
W_{d,3}\Theta_{\cO}=\Theta_{\cO},\qquad W_{4d,3}\Theta_{\cO_2}=\Theta_{\cO_2},\qquad
W_{4d,3}\Theta_\cO=8\,\Theta_\cO(4\,\cdot),
\]
and similarly $W_{4d,1}\vartheta_{\cO_2}=\vartheta_{\cO_2}$ and
$W_{4d,1}\vartheta_\cO=2\,\vartheta_\cO(4\,\cdot)$. Define
\begin{equation}\label{eq:F2F0}
F_2=3\Theta_{\cO_2}-\Theta_\cO-8\Theta_\cO(4z),\qquad
F_0=3\vartheta_{\cO_2}-\vartheta_\cO-2\vartheta_\cO(4z).
\end{equation}
Then $W_{4d,3}F_2=F_2$, $W_{4d,1}F_0=F_0$, and the constant term of $F_0$ is
$\frac32-\frac12-1=0$.

\begin{lemma}\label{lem:traces}
\[
\sum_{\eta\in\widehat{\Cl}}g_{\psi\eta}=h\,\Theta_\cO,\qquad
\sum_{\chi\in\cR}g_{\psi\chi}=h\,F_2,\qquad
\sum_{\rho\in\cP}h_\rho=\frac h2F_0 .
\]
\end{lemma}

\begin{proof}
The first identity is character orthogonality on $\Cl(K)$. For the second, summing the
theta series restricted to ideals coprime to $2$ over all $3h$ ray characters produces
$3h\,\Theta_{1+2\cO}$. For $\eta\in\widehat{\Cl}$, the restricted series is
$g_{\psi\eta}(z)-4g_{\psi\eta}(4z)$ because $\psi\eta((2))=4$. Substituting
$\Theta_{1+2\cO}=\Theta_{\cO_2}-4\Theta_\cO(4z)$ gives
\[
\sum_{\chi\in\cR}g_{\psi\chi}=3h\bigl(\Theta_{\cO_2}-4\Theta_\cO(4z)\bigr)-h\bigl(\Theta_\cO-4\Theta_\cO(4z)\bigr)=hF_2 .
\]
The weight-one calculation is identical, with $\eta((2))=1$ and
$\vartheta_{2\cO}=\vartheta_\cO(4z)$, giving $\sum_{\chi\in\cR}h_\chi=hF_0$; applying
\eqref{eq:conj} divides this sum by $2$.
\end{proof}

\begin{proposition}\label{prop:FE}
\begin{align*}
\sum_{\eta\in\widehat{\Cl}}L'(g_{\psi\eta},0)&=\frac{d^{3/2}}{4\pi^3}M_{\widehat{\Cl}},&
\sum_{\chi\in\cR}L'(g_{\psi\chi},0)&=\frac{(4d)^{3/2}}{4\pi^3}M_{\cR},\\
\sum_{\rho\in\cP}L'(h_\rho,-1)&=-\frac{d^{3/2}}{\pi^3}Z_{\cP}.
\end{align*}
\end{proposition}

\begin{proof}
The invariance of the complete traces under $W_{N,k}$ gives their functional equations,
and comparing the two relevant endpoints of each Mellin transform determines the constants.

Let $f$ be a cusp form with $W_{N,k}f=f$ and real Fourier coefficients. Its completed
Mellin transform $\Lambda_N(s)=N^{s/2}\int_0^\infty f(iy)y^{s-1}dy=N^{s/2}(2\pi)^{-s}\Gamma(s)L(f,s)$
is entire and satisfies $\Lambda_N(s)=\Lambda_N(k-s)$.

When $k=3$, the simple pole of $\Gamma$ at $0$ gives $\Lambda_N(0)=L'(f,0)$, whereas
$\Lambda_N(3)=N^{3/2}L(f,3)/(4\pi^3)$. When $k=1$, the residue of $\Gamma$ at $-1$ is
$-1$, so $\Lambda_N(-1)=-2\pi N^{-1/2}L'(f,-1)$, while $\Lambda_N(2)=NL(f,2)/(4\pi^2)$;
hence $L'(f,-1)=-N^{3/2}L(f,2)/(8\pi^3)$. Applying these identities to $h\Theta_\cO$
($N=d$), $hF_2$ ($N=4d$), and $\frac h2F_0$ ($N=4d$) via Lemma~\ref{lem:traces} proves the
proposition, with exponential decay at both endpoints guaranteed by the modular
transformation laws.
\end{proof}

\begin{theorem}\label{thm:GT}
Theorem~A holds:
\[
\begin{aligned}
n_4(s_d)={}&\frac{40}{dh}\sum_{\eta\in\widehat{\Cl}}L'(g_{\psi\eta},0)
+\frac{40}{3dh}\sum_{\chi\in\cR}L'(g_{\psi\chi},0)\\[4pt]
&+\frac{5\sqrt d}{h\pi^3}\sum_{\eta\in\widehat{\Cl}}L_K(\eta,2)
-\frac{80}{3dh}\sum_{\rho\in\cP}L'(h_\rho,-1).
\end{aligned}
\]
\end{theorem}

\begin{proof}
By Lemma~\ref{lem:domain}, $\tau_d\in\Cinf$, so \eqref{eq:lattice-d} applies.
Substituting Proposition~\ref{prop:decomp} and then Proposition~\ref{prop:FE}, the three
prefactors simplify to
\[
\frac{5\sqrt d}{h\pi^3}\cdot2\cdot\frac{4\pi^3}{d^{3/2}}=\frac{40}{dh},\qquad
\frac{5\sqrt d}{h\pi^3}\cdot\frac{16}3\cdot\frac{4\pi^3}{(4d)^{3/2}}=\frac{40}{3dh},\qquad
\frac{5\sqrt d}{h\pi^3}\cdot\frac{16}3\cdot\Bigl(-\frac{\pi^3}{d^{3/2}}\Bigr)=-\frac{80}{3dh}.
\qedhere
\]
\end{proof}

The proof makes no assumption of class number one, no splitting of $\Cl_{(2)}(K)$ as
$C_3\times\Cl(K)$, and no requirement that the characters in $\cR$ have order three.

\Needspace{210pt}
\subsection{Individual root numbers}
\begin{proposition}\label{prop:signs}
Let $f$ be one of the newforms $g_{\psi\eta}$ ($\eta\in\widehat{\Cl}$),
$g_{\psi\chi}$ ($\chi\in\cR$), or $h_\rho$ ($\rho\in\cP$). Let $k$ be its weight
($3$ or $1$), $N$ its level, and $\varepsilon_f$ its sign under $W_{N,k}$.
Then $\varepsilon_f=+1$.
\end{proposition}

\begin{proof}
For $\Psi=\psi\eta$ or $\psi\chi$, $|\Psi(\mathfrak a)|=\Nm\mathfrak a$ and
$\Psi(\mathfrak a)\Psi(\bar{\mathfrak a})=(\Nm\mathfrak a)^2$, so
$\Psi(\bar{\mathfrak a})=\overline{\Psi(\mathfrak a)}$. Because $\Nm\bar{\mathfrak a}=\Nm\mathfrak a$,
the Fourier coefficients of $g_\Psi$ are real. For the finite characters,
$\rho(\bar{\mathfrak a})=\rho(\mathfrak a)^{-1}=\overline{\rho(\mathfrak a)}$ likewise
implies that $h_\rho$ has real coefficients. By \cite[Theorem~1.1]{AtkinLi78}, the
Fricke operator sends a normalized primitive newform to a scalar multiple of the
newform with complex-conjugate coefficients, which here is $f$ itself. Hence $f$ is an
eigenform of the involution $W_{N,k}$, with eigenvalue $\varepsilon_f=\pm1$.

Each of the fixed traces $h\Theta_\cO$, $hF_2$, and $\frac h2F_0$ decomposes as a sum
$\sum_jc_jf_j$ with $c_j>0$ and normalized first coefficient $a_1(f_j)=1$. Applying
$W_{N,k}-1$ gives $\sum_jc_j(\varepsilon_j-1)f_j=0$, and reading off the coefficient of $q$
yields $\sum_jc_j(1-\varepsilon_j)=0$. Since every summand is nonnegative, each
$\varepsilon_j$ must equal $1$.
\end{proof}

\begin{remark}
Because $W_{N,k}$ differs from the standard Fricke operator by a factor of $i^k$, having
$\varepsilon_f=+1$ in odd weight does not mean that the unnormalized Fricke operator has
eigenvalue $+1$.
\end{remark}

\section{Class number one: the cubic parameter fields}\label{sec:classone}

In this section let $d\in S_0=\{11,19,43,67,163\}$, so that $h(K)=1$ and
$\Cl_{(2)}(K)\cong(\cO/2\cO)^\times\cong C_3$. Fix a nontrivial character $\chi_d$
normalized by $\chi_d((\alpha))=\chi_d^{\rm res}(\bar\alpha)$ with
$\chi_d^{\rm res}(\bar\tau_d)=e^{2\pi i/3}$. Then $\cR=\{\chi_d,\chi_d^{-1}\}$ and
$\cP=\{\chi_d\}$, and by \eqref{eq:conj} there is a single weight-one newform
\[
h_{4d}:=h_{\chi_d}=h_{\chi_d^{-1}}=\sum_{(\mathfrak a,2)=1}\chi_d(\mathfrak a)q^{\Nm\mathfrak a}.
\]
Set $H_{4d}=L'(h_{4d},-1)$.

\subsection{The cubic field}

\begin{proposition}\label{prop:cubic}
Let $F_d=\Q(s_d)$. Then $F_d$ is a cubic field of signature $(1,1)$ and discriminant
$-4d$, and
\begin{equation}\label{eq:zetaF}
\zeta_{F_d}(s)=\zeta(s)\,L(h_{4d},s).
\end{equation}
Consequently, $\zeta'_{F_d}(-1)=-\tfrac1{12}H_{4d}$.
\end{proposition}

\begin{proof}
The proof identifies a small defining cubic for $\Q(s_d)$, computes its ring of
integers, and matches its Artin $L$-function with the cubic ray character $\chi_d$. Differentiating at
$s=-1$ then gives the derivative identity.

\smallskip
\Needspace{90pt}
\noindent\textbf{Field identification.}\par
Let $p_d$ and $r_d(x)=u_dx^2+v_dx+w_d$ be as in Table~\ref{tab:fields}. Exact
polynomial division shows that $Q_{j_d}(r_d(x))\equiv0\pmod{p_d(x)}$. Each $p_d$ is
irreducible by the rational root test and has negative discriminant; let $\alpha_d$ be
its unique real root. Since $u_d\ne0$, $r_d(\alpha_d)$ is irrational and therefore
generates the cubic field $\Q(\alpha_d)$. As this field has signature $(1,1)$,
$r_d(\alpha_d)$ must be the unique real root of $Q_{j_d}$, which is $s_d$ by
Proposition~\ref{prop:minpoly}. Thus $F_d=\Q(\alpha_d)$.

\smallskip
\Needspace{90pt}
\noindent\textbf{Discriminant.}\par
For $d=163$, the algebraic integer $v=(\alpha^2-2\alpha+5)/2$ satisfies $v^3-8v-10=0$,
whose discriminant is $-652$. For the remaining values of $d$, the polynomial $p_d$
already has discriminant $-4d$. In every case, the only possible index primes are $2$
(together with $3$ for $d=27$, needed in Section~\ref{sec:d27}). Dedekind's criterion
rules them out via the congruence $P_d\equiv(x-r)^3\pmod\ell$ with
$P_d(r)/\ell\not\equiv0\pmod\ell$, witnessed by $(\ell,r)=(2,1)$ for $d=11,43,67$ and
$(2,0)$ for $d=19,163$ (each with $P_d(r)/\ell\equiv1$), and by $(2,0)$ and $(3,-1)$
for $d=27$. Hence $D_{F_d}=-4d$.

\smallskip
\Needspace{90pt}
\noindent\textbf{Artin decomposition.}\par
Let $L_d$ be the Galois closure of $F_d$. Because $D_{F_d}$ is negative and not a
square, $\Gal(L_d/\Q)\cong S_3$ with quadratic subfield $\Q(\sqrt{-4d})=K$. The
permutation representation of $F_d$ decomposes as $1\oplus\rho$, so
$\zeta_{F_d}(s)=\zeta(s)L(\rho,s)$.

The extension $L_d/K$ is cyclic of degree three. Its two nontrivial characters share a
conductor $\mathfrak f$, and the conductor--discriminant formula gives
$|D_{L_d}|=|D_K|^3\Nm(\mathfrak f)^2$. On the other hand, the regular representation
of $S_3$ is $1\oplus\varepsilon\oplus2\rho$, so multiplicativity of conductors
\cite[Section~2.5]{JonesRoberts17} yields $|D_{L_d}|=|D_K|\,|D_{F_d}|^2$. Thus
$\Nm\mathfrak f=|D_{F_d}|/|D_K|=4$, which forces $\mathfrak f=(2)$ since $2$ is inert.
Because $\Cl_{(2)}(K)$ has order three, the characters of $L_d/K$ are $\chi_d^{\pm1}$,
and induction \cite[Ch.~VIII, Section~10]{MilneCFT} gives
$L(\rho,s)=L_K(\chi_d,s)=L(h_{4d},s)$.

The ramified Euler factors also coincide. At $2$, inertia in $L_d/\Q$ contains $C_3$,
so $\rho^{I_2}=0$, matching the Hecke factor $1$ at the ramified prime $(2)$. At $p=d$,
inertia is a subgroup $C_2$ equal to the decomposition group, and Frobenius acts
trivially on the one-dimensional invariant subspace, giving $(1-p^{-s})^{-1}$; on the
Hecke side, the prime above $d$ is $(\sqrt{-d})$ and $\sqrt{-d}=2\tau_d-1\equiv1\pmod2$,
producing the same factor.

Finally, $L(h_{4d},-1)=0$ and $\zeta(-1)=-1/12$, so differentiating \eqref{eq:zetaF}
at $s=-1$ yields $\zeta'_{F_d}(-1)=-\tfrac1{12}H_{4d}$.
\end{proof}

\begin{table}[ht]
\centering
\small
\setlength{\tabcolsep}{8pt}
\renewcommand{\arraystretch}{1.15}
\begin{tabular}{@{}rlr@{}}
\toprule
$d$ & $p_d(x)$ & $D_{F_d}$\\
\midrule
$11$ & $x^3-x^2+x+1$ & $-44$\\
$19$ & $x^3-2x-2$ & $-76$\\
$27$ & $x^3-2$ & $-108$\\
$43$ & $x^3-x^2-x+3$ & $-172$\\
$67$ & $x^3-x^2-3x+5$ & $-268$\\
$163$ & $x^3-x^2+7x+5$ & $-652$\\
\bottomrule
\end{tabular}

\medskip
\begin{tabular}{@{}rrrr@{}}
\toprule
$d$ & $u_d$ & $v_d$ & $w_d$\\
\midrule
$11$ & $-11\,328$ & $17\,344$ & $-20\,624$\\
$19$ & $-119\,232$ & $-212\,544$ & $-136\,080$\\
$27$ & $-2\,579\,200$ & $-3\,252\,800$ & $-4\,096\,144$\\
$43$ & $-121\,824\,000$ & $287\,452\,800$ & $-268\,905\,744$\\
$67$ & $-12\,375\,633\,600$ & $36\,131\,832\,000$ & $-32\,233\,449\,744$\\
$163$ & $-27\,882\,300\,234\,225\,600$ & $45\,279\,829\,361\,980\,800$ & $-223\,429\,048\,349\,227\,344$\\
\bottomrule
\end{tabular}
\caption{Cubic fields and parameter coordinates, displayed in separate panels.
$s_d=u_d\alpha_d^2+v_d\alpha_d+w_d$, where $p_d(\alpha_d)=0$;
for $d=163$, $\operatorname{disc}p_d=-2608$.}
\label{tab:fields}
\end{table}

\begin{remark}[The lattice form of the ray-class constant]\label{rem:Dd}
Set $c_d=(d+1)/4$ and assign to $(a,b)\in\Z^2$ the parity weight $w(a,b)=0$ if
$a\equiv b\equiv0\pmod2$, $w(a,b)=2$ if $a$ is odd and $b$ is even, and $w(a,b)=-1$
if $b$ is odd. Partitioning the lattice $\cO$ gives the $q$-series identity
\begin{equation}\label{eq:h4d-theta}
h_{4d}(z)=\frac14\sum_{a,b}w(a,b)q^{a^2+ab+c_db^2}
=\frac12\Bigl(\sum_{u,v}q^{u^2+dv^2}-\sum_{u,v}q^{4u^2+2uv+c_dv^2}\Bigr).
\end{equation}
Indeed, if $T$ is the theta series of $\cO$, $E$ is that of $2\cO$, and $A,B,C$ are the
subsums where $b$, $a$, and $a+b$ are even, then $A+B+C=T+2E$ and conjugation gives
$B=C$, so the half-weighted sum
$W=\frac12\sum_{a,b}w(a,b)q^{a^2+ab+c_db^2}$ satisfies
$W=\frac12(3A-2E-T)=A-B$; substituting $b=2v$, $u=a+v$ in
$A$ and $a=2u$ in $B$ gives the right-hand side. On the other hand, every ideal coprime
to $2$ has two generators and $\chi_d+\chi_d^{-1}$ takes the values $2,-1,-1$ on the
three nonzero residues modulo $2\cO$, so $W=2h_{4d}$.

Consequently, the ray-class constant
\[
\cD_d=\frac{3d^{3/2}}{2\pi^3}\cdot\frac12\sideset{}{'}\sum w(a,b)(a^2+ab+c_db^2)^{-2}
\]
satisfies
\[
\cD_d=\frac{3d^{3/2}}{\pi^3}L(h_{4d},2)=-3H_{4d}=36\,\zeta'_{F_d}(-1),
\]
using $H_{4d}=-(4d)^{3/2}L(h_{4d},2)/(8\pi^3)$ from Proposition~\ref{prop:FE}. This
$\cD_d$ is precisely the new constant appearing in \cite[Section~9]{Zheng26}. Since
$L(h_{4d},s)=\zeta_{F_d}(s)/\zeta(s)$ is attached to an irreducible two-dimensional
representation, it is not expected to reduce to products of Dirichlet $L$-values.
\end{remark}

\subsection{The five fundamental evaluations}

\begin{theorem}\label{thm:six}
For $d\in S_0$, with $M_d,U_d,\Delta_d$ as in \eqref{eq:consts},
\[
n_4(s_d)=\frac{10}{9d}\Bigl(36M_d+12U_d+3\Delta_d+288\,\zeta'_{F_d}(-1)\Bigr).
\]
\end{theorem}

\begin{proof}
Specializing Theorem~\ref{thm:GT} to $h=1$ gives
\[
n_4(s_d)=\frac{40}dM_d+\frac{40}{3d}U_d+\frac{5\sqrt d}{\pi^3}\zeta_K(2)-\frac{80}{3d}H_{4d}.
\]
The completed $L$-function $(d/\pi)^{(s+1)/2}\Gamma(\frac{s+1}2)L(\chi_{-d},s)$ of the
odd character $\chi_{-d}$ is invariant under $s\mapsto1-s$; comparing $s=-1$ and $s=2$
yields $\Delta_d=d^{3/2}L(\chi_{-d},2)/(4\pi)$. Combining this with
$\zeta_K(2)=\frac{\pi^2}6L(\chi_{-d},2)$ turns the third term into
\[
\frac{5\sqrt d}{\pi^3}\cdot\frac{\pi^2}6\cdot\frac{4\pi}{d^{3/2}}\Delta_d=\frac{10}{3d}\Delta_d .
\]
Proposition~\ref{prop:cubic} gives $-\frac{80}{3d}H_{4d}=\frac{320}d\zeta'_{F_d}(-1)$,
and factoring out $10/(9d)$ yields the stated formula.
\end{proof}

Here $g_\psi$ is the classical weight-three CM newform of level $d$ attached to
$\psi((\alpha))=\alpha^2$, while $g_{\psi\chi_d^{\pm1}}$ are the two newforms of level
$4d$; by Proposition~\ref{prop:signs}, all three have sign $+1$.

\section{The non-fundamental discriminant \texorpdfstring{$-27$}{-27}}\label{sec:d27}

For $d=27$ $\tau_{27}=(1+3\sqrt{-3})/2=3\omega-1$ with $\omega=(1+\sqrt{-3})/2$,
so the two lattices in \eqref{eq:lattice} are $R_3$ and $R_6$, where $R_f=\Z+f\omega\Z$.
Put $K=\Q(\sqrt{-3})$, $\cO=\Z[\omega]$ (so $\omega^2=\omega-1$ and $\cO^\times=\mu_6$),
and $\mathfrak p=(2\omega-1)=(\sqrt{-3})$, so $(3)=\mathfrak p^2$. Because $-27$ is
non-fundamental, Theorem~A does not apply directly; the key step is to separate the
lattice theta series into its primitive constituents.

\subsection{Primitive characters}
Reduction modulo $3$ gives an isomorphism $\mu_6\cong(\cO/3\cO)^\times$; for $\alpha$
coprime to $3$, let $u_3(\alpha)\in\mu_6$ be the unique unit congruent to $\alpha$
modulo $3$. The group $(\cO/2\cO)^\times$ is cyclic of order three. Define
$\chi_\pm(\bar\omega)=\omega^{\pm2}$ and set
\[
\varepsilon_{12}(\alpha)=\chi_-(\bar\alpha),\qquad
\varepsilon_{27}(\alpha)=u_3(\alpha)^{-2},\qquad
\varepsilon_{108}(\alpha)=u_3(\alpha)^2\chi_+(\bar\alpha),
\]
for $\alpha$ coprime to $2$, $3$, and $6$, respectively.

\begin{lemma}\label{lem:psiN}
The formulas $\psi_N((\alpha))=\alpha^2\varepsilon_N(\alpha)$, $N\in\{12,27,108\}$, define
primitive Hecke characters of infinity type $2$ and conductors $(2)$, $\mathfrak p^2=(3)$,
and $(6)$, respectively. Their theta series $g_N$ are normalized newforms of weight $3$,
level $N$, and nebentypus $\chi_{-3}$.
\end{lemma}

\begin{proof}
For any unit $u\in\mu_6$, $\chi_\pm(u)=u^{\pm2}$ and $u_3(u)=u$, so
$\varepsilon_N(u)=u^{-2}$. Thus $\alpha^2\varepsilon_N(\alpha)$ depends only on the ideal
$(\alpha)$, and multiplicativity is immediate.

For $\psi_{12}$, nontriviality on $\mathbb F_4^\times$ excludes conductor $1$. For
$\psi_{27}$, the definition of $u_3$ implies that its conductor divides $(3)=\mathfrak p^2$;
since $\omega^2-1=(2\omega-1)\omega$, $\omega^2\equiv1\pmod{\mathfrak p}$ while
$\varepsilon_{27}(\omega^2)=\omega^2\ne1$, so the conductor is precisely $\mathfrak p^2$
of norm $9$, giving level $3\cdot9=27$. For $\psi_{108}$, the elements
$\beta_3=3+2\omega$ and $\beta_2=1+3\omega$ (of norms $19$ and $13$) satisfy
\[
\beta_3\equiv1\pmod{2\mathfrak p},\ \varepsilon_{108}(\beta_3)=\omega^2\ne1;\qquad
\beta_2\equiv1\pmod 3,\ \varepsilon_{108}(\beta_2)=\omega^4\ne1,
\]
excluding the two maximal proper divisors of $(6)$.

Hecke's theorem \cite[Theorem~3.7]{Fei21} therefore gives weight $3$ and levels
$|D_K|\Nm\mathfrak f\in\{12,27,108\}$. Since all three finite parts equal $1$ on
rational integers coprime to their conductors, the nebentypus is $\chi_{-3}$.
\end{proof}

The cubic ideal character
\begin{equation}\label{eq:phi}
\phi((\alpha))=u_3(\alpha)^{-2}\chi_+(\bar\alpha),\qquad(\alpha,6)=1,
\end{equation}
likewise has conductor $(6)$ by the same witnesses, and links the three primitive Hecke
characters on ideals coprime to $6$ by cubic twists:
\[
\psi_{12}\ \bigl(\text{cond.\ }(2)\bigr)
\xleftarrow{\ \cdot\,\phi^{-1}\ }
\psi_{27}\ \bigl(\text{cond.\ }\mathfrak p^2=(3)\bigr)
\xrightarrow{\ \cdot\,\phi\ }
\psi_{108}\ \bigl(\text{cond.\ }(6)\bigr).
\]

\subsection{The lattice theta series}
Assign to the residues modulo $2$ the weights $w(0)=0$, $w(1)=2$, and
$w(\bar\omega)=w(\bar\omega^2)=-1$, and put
\begin{align*}
F^{\rm lat}(z)&=\frac12\sum_{\alpha\in R_3}
 w(\alpha)\alpha^2q^{\Nm\alpha},\\[4pt]
F^{(6)}(z)&=\frac12\sum_{\substack{\alpha\in R_3\\(\alpha,6)=1}}
 w(\alpha)\alpha^2q^{\Nm\alpha}.
\end{align*}
Here $(\alpha,6)=1$ means that $\alpha\cO$ is coprime to $6\cO$; since $w(0)=0$ already
eliminates terms divisible by $2$, the second sum is equivalently restricted to
$3\nmid\Nm\alpha$. Write $\Theta_f=\Theta_{R_f}$.

\begin{proposition}\label{prop:F-lat}
As holomorphic functions,
\begin{align}
g_{27}&=\Theta_3,\qquad g_{12}=\Theta_2=\eta(2z)^3\eta(6z)^3,\label{eq:g27g12}\\
F^{(6)}&=g_{108}+g_{12}+3g_{12}(3z),\label{eq:F6}\\
F^{\rm lat}&=g_{108}+g_{12}+3g_{12}(3z)+27g_{12}(9z).\label{eq:Flat}
\end{align}
Put $C_N=N^{3/2}/(4\pi^3)$ and $M_N=L'(g_N,0)$. Then $M_N=C_NL(g_N,3)$ and
\begin{equation}\label{eq:Flat-L}
L'(F^{\rm lat},0)=C_{108}L(F^{\rm lat},3)=M_{108}+31M_{12},
\end{equation}
whereas $C_{108}L(F^{(6)},3)=M_{108}+30M_{12}$ and $L'(F^{(6)},0)=M_{108}+4M_{12}$.
\end{proposition}

\begin{proof}
First identify the primitive newforms $g_{27}$ and $g_{12}$, then isolate the terms
coprime to $6$, and finally restore the terms divisible by $3$, which supply the
oldform $27g_{12}(9z)$ responsible for the final coefficient $31$.

\smallskip
\Needspace{90pt}
\noindent\textbf{The primitive forms $g_{27}$ and $g_{12}$.}\par
Among the six generators of an ideal coprime to $3$, exactly two are congruent to
$\pm1\pmod3$, and these are its generators in $R_3$. Writing $\alpha=a+3b\omega$,
$3\mid\Nm\alpha$ if and only if $3\mid a$, so the remaining elements form $3\cO$; their
harmonic theta series vanishes because $\sum_{u\in\mu_6}u^2=0$. Hence $g_{27}=\Theta_3$,
and the analogous argument modulo $2$ gives $g_{12}=\Theta_2$. Since
$\eta(2z)^3\eta(6z)^3\in S_3(\Gamma_0(12),\chi_{-3})$, checking coefficients up to the
Sturm bound $6$ \cite{Sturm87} establishes the eta-product identity, already recorded
in \cite[(1.2)]{Samart13}.

\smallskip
\Needspace{90pt}
\noindent\textbf{The identity \eqref{eq:F6}.}\par
On a primary generator modulo $3$, $\chi_+(\alpha)+\chi_-(\alpha)=w(\alpha)$, so
$F^{(6)}$ equals the sum of $g_{108}$ and the theta series of $\psi_{12}$ with the Euler
factor at $\mathfrak p$ removed. Since $\psi_{12}(\mathfrak p)=(2\omega-1)^2=-3$, this
gives \eqref{eq:F6}.

\smallskip
\Needspace{90pt}
\noindent\textbf{The identity \eqref{eq:Flat}.}\par
The terms omitted from $F^{\rm lat}$ in $F^{(6)}$ are those of the form $3\beta$ with
$\beta\in\cO$. Their residue sums modulo $2\cO$ satisfy $S_1+S_a+S_b=0$, whence
$2S_1-S_a-S_b=3S_1=6g_{12}(z)$; multiplying by $(3\beta)^2q^{9\Nm\beta}$ contributes
$27g_{12}(9z)$.

\smallskip
\Needspace{90pt}
\noindent\textbf{Functional equations and constants.}\par
By Lemma~\ref{lem:poisson}, $g_{12}$, $g_{27}$, and $\Theta_6$ are invariant under
$W_{12,3}$, $W_{27,3}$, and $W_{108,3}$, respectively. Partitioning modulo $2$ gives
$F^{\rm lat}=3\Theta_6-g_{27}-8g_{27}(4z)$, and since $W_{108,3}$ exchanges $g_{27}$
and $8g_{27}(4z)$, it fixes $F^{\rm lat}$. For each $r\mid9$,
$W_{108,3}\{g_{12}(r\,\cdot)\}(iy)=\frac{27}{r^3}g_{12}(\frac9riy)$, so
$g_{12}+3g_{12}(3z)+27g_{12}(9z)$ is also fixed, which implies $W_{108,3}g_{108}=g_{108}$.
Proposition~\ref{prop:FE} then gives $M_N=C_NL(g_N,3)$. Passing to $L$-series from
\eqref{eq:F6} and \eqref{eq:Flat} yields
\begin{align*}
L(F^{\rm lat},s)&=L(g_{108},s)+(1+3^{1-s}+3^{3-2s})L(g_{12},s),\\
L(F^{(6)},s)&=L(g_{108},s)+(1+3^{1-s})L(g_{12},s),
\end{align*}
and using $C_{108}/C_{12}=27$ together with $L(g_{12},0)=0$ gives the stated constants.
\end{proof}

Notice that the two values attached to $F^{(6)}$ differ by $26M_{12}>0$ and must not
be conflated. Algebraically, the Dirichlet polynomial $P(s)=1+3^{1-s}+3^{3-2s}$
satisfies $P(s)=27\cdot9^{-s}P(3-s)$, which explains why the full lattice series obeys
the level-$108$ functional equation while the truncated series does not.

\subsection{The scalar term and the evaluation}
Let $X_{27}=A(\omega+R_3)$. The number of generators of an ideal whose $\omega$-coefficient
is congruent to $1\pmod3$ is $2$, $3$, or $0$ according as its $\mathfrak p$-valuation is
$0$, $1$, or at least $2$. Therefore,
\begin{equation}\label{eq:X27}
X_{27}=(2+3\cdot3^{-2})(1-3^{-2})\zeta_K(2)=\frac{56}{27}\zeta_K(2),\qquad
A(R_3)=6\zeta_K(2)-2X_{27}=\frac{50}{27}\zeta_K(2).
\end{equation}
In particular, the constant $9\sqrt3X_{27}/(2\pi^3)$ occurring in
\cite[Remark~9.3]{Zheng26} equals $\frac{56}{27}\Delta_3$ exactly, where
$\Delta_3=L'(\chi_{-3},-1)=\frac{3\sqrt3}{4\pi}L(\chi_{-3},2)$.
This proves the simplification observed numerically there.

Let $A_r,H_r$ ($r=1,a,b$) denote the scalar and harmonic sums over the three nonzero
cosets of $2R_3$ in $R_3$, and set
\[
BU_{27}=\tfrac12(2A_1-A_a-A_b),\qquad U^{\rm lat}_{27}=\tfrac12(2H_1-H_a-H_b),\qquad
\cD_{27}=\frac{9\sqrt3}{2\pi^3}BU_{27}.
\]
Then $U^{\rm lat}_{27}=L(F^{\rm lat},3)$, and since $R_6=2R_3\sqcup(1+2R_3)$, this gives
\begin{equation}\label{eq:T27}
-T(R_3)+4T(R_6)=\tfrac12A(R_3)+H(R_3)+\tfrac83BU_{27}+\tfrac{16}3U^{\rm lat}_{27},
\qquad H(R_3)=2L(g_{27},3).
\end{equation}

\begin{proposition}\label{prop:cubic27}
$F_{27}=\Q(s_{27})\cong\Q(\sqrt[3]2)$, with
\[
s_{27}=-4\,096\,144-3\,252\,800\sqrt[3]2-2\,579\,200\sqrt[3]4,
\]
and $\cD_{27}=\frac43\zeta'_{\Q(\sqrt[3]2)}(-1)$.
\end{proposition}

\begin{proof}
The field identification is row $d=27$ of Table~\ref{tab:fields}, and the maximality
witnesses at $2$ and $3$ in Proposition~\ref{prop:cubic} give $D_{F_{27}}=-108$.

The Galois closure of $F_{27}$ has quadratic subfield $K$, and the cubic characters of
its relative extension have conductor of norm $108/3=36$, which forces the conductor to
be $(6)$. Since $(\cO/6\cO)^\times$ has order $18$ and the image of $\mu_6$ has order $6$,
the ray class group modulo $(6)$ has order $3$, so these characters are $\phi^{\pm1}$.
Inertia at both $2$ and $3$ contains $C_3$, so the Euler factors at $2$ and $3$ equal $1$
on both sides, giving $\zeta_{F_{27}}(s)=\zeta(s)L(h_{108},s)$ with
$h_{108}=\sum_{(\mathfrak a,6)=1}\phi(\mathfrak a)q^{\Nm\mathfrak a}$.

Ideals coprime to $6$ have two generators in $R_3$ and contribute $\phi+\phi^{-1}$ to
the weighted sum, while the remaining elements with $3\mid\Nm\alpha$ form $3\cO$, whose
weighted scalar contribution vanishes since $2-1-1=0$. Thus $BU_{27}=2L(h_{108},2)$.
The completed functional equation of $L(h_{108},s)=\zeta_{F_{27}}(s)/\zeta(s)$
has sign $+1$ and conductor $108$. The weight-one calculation in the proof of
Proposition~\ref{prop:FE} therefore gives
$H_{108}=L'(h_{108},-1)=-108^{3/2}L(h_{108},2)/(8\pi^3)$, whence
$\cD_{27}=-\tfrac19H_{108}=\tfrac43\zeta'_{F_{27}}(-1)$. (For reference,
$h_{108}=\eta(6z)\eta(18z)$ by \cite[Lemma~2.2]{HLMYYZ23}.)
\end{proof}

\begin{theorem}\label{thm:d27}
Theorem~C holds:
\[
n_4(s_{27})=\frac{10}{81}\Bigl(12M_{27}+4M_{108}+124M_{12}+25\Delta_3+96\,\zeta'_{\Q(\sqrt[3]2)}(-1)\Bigr).
\]
\end{theorem}

\begin{proof}
Since $\Ima\tau_{27}=3\sqrt3/2$, Proposition~\ref{prop:lattice} gives
$n_4(s_{27})=\frac{15\sqrt3}{\pi^3}(-T(R_3)+4T(R_6))$. Substituting \eqref{eq:T27},
\eqref{eq:X27}, and $M_{27}=C_{27}L(g_{27},3)$ yields the coefficients
\[
\frac{15\sqrt3}{\pi^3}\cdot2\cdot\frac{4\pi^3}{81\sqrt3}=\frac{40}{27},\qquad
\frac{15\sqrt3}{\pi^3}\cdot\frac{16}3\cdot\frac{\pi^3}{162\sqrt3}=\frac{40}{81},\qquad
\frac{15\sqrt3}{\pi^3}\cdot\frac{25}{27}\cdot\frac{2\pi^3}{9\sqrt3}=\frac{250}{81},
\]
for $M_{27}$, $L'(F^{\rm lat},0)$, and $\Delta_3$, together with
$\frac{15\sqrt3}{\pi^3}\cdot\frac83\cdot\frac{2\pi^3}{9\sqrt3}=\frac{80}9$ for $\cD_{27}$.
Inserting \eqref{eq:Flat-L} and Proposition~\ref{prop:cubic27} completes the proof.
\end{proof}

\begin{remark}
With the notation of \cite{Zheng26}, the lattice normalization is
\[
MU_{27}=C_{108}L(F^{\rm lat},3)=M_{108}+31M_{12}.
\]
Thus the correction is $31M_{12}$, rather than $30M_{12}$ or $4M_{12}$. This
clarification does not assert an error in \cite{Zheng26}.
\end{remark}

\section{The first case of class number two: \texorpdfstring{$d=35$}{d=35}}\label{sec:d35}

Let $K=\Q(\sqrt{-35})$, so that $h=2$ and $\Cl(K)=\{[\cO],[\mathfrak p_3]\}$, represented
by the reduced quadratic forms $[1,1,9]$ and $[3,1,3]$. The nontrivial class character
is the genus character $\eta_5$ associated with $\Q(\sqrt5,\sqrt{-7})$, giving
\begin{equation}\label{eq:genus35}
L_K(\eta_5,s)=L(s,\chi_5)L(s,\chi_{-7}),\qquad \zeta_K(s)=\zeta(s)L(s,\chi_{-35}).
\end{equation}

\subsection{The ray class field modulo \texorpdfstring{$2$}{2}}

\begin{lemma}\label{lem:ray35}
The ray class group $\Cl_{(2)}(K)$ is cyclic of order six and coincides with the
Picard group of $\Z[\sqrt{-35}]$. Its class field $N$ has degree $12$ over $\Q$ with
$\Gal(N/\Q)\cong D_6$, and
\[
N=C(\sqrt5,\sqrt{-7}),\qquad C=\Q(a),\quad a^3+2a-2=0 .
\]
The field $C$ has discriminant $-140$ and one complex place, while
$E:=C(\sqrt5)$ has discriminant $98\,000$ and signature $(2,2)$.
\end{lemma}

\begin{proof}
Because $2$ is inert, $(\cO/2\cO)^\times\cong\mathbb F_4^\times$ has order $3$ and
$\cO^\times=\{\pm1\}$ maps to $1$, so $|\Cl_{(2)}(K)|=3h=6$ (Section~\ref{ssec:ray}). An
element coprime to $2$ is congruent to $1\pmod{2\cO}$ if and only if it is congruent to
a rational integer modulo $2\cO$, since every odd integer is $\equiv1\pmod2$. Thus
$\Cl_{(2)}(K)$ is the ring class group $\mathrm{Pic}(\Z[\sqrt{-35}])$ of conductor $2$.
The reduced forms of discriminant $-140$ are
\[
[1,0,35],\quad[5,0,7],\quad[3,\pm2,12],\quad[4,\pm2,9].
\]
Only one non-principal form is ambiguous, so $\Cl_{(2)}(K)$ has a unique element of
order two and is cyclic. The forms $[4,\pm2,9]$ represent $4$, which is a square modulo
$5$ and $7$; they lie in the principal genus and have order three, whereas
$[3,\pm2,12]$ have order six.

As a ring class field, $N$ is generalized dihedral over $\Q$, so $\Gal(N/\Q)\cong D_6$,
and it contains the genus field $\Q(\sqrt5,\sqrt{-7})$. Its cubic subfields are the
cubic fields of discriminant $-140$. Since $x^3+2x-2$ is Eisenstein at $2$ and has
discriminant $-4\cdot8-27\cdot4=-140$, its field discriminant is $-140$ ($-35$ cannot
be a cubic field discriminant because $h(-35)=2$). Thus $C\subset N$ up to conjugacy,
and comparing degrees gives $N=C(\sqrt5,\sqrt{-7})$. The conductor--discriminant formula
and Lemma~\ref{lem:Artin35} yield $|D_E|=|D_C|\cdot5\cdot140=98\,000$.
\end{proof}

Fix a generator $\psi_3$ of the order-three characters of $\Cl_{(2)}(K)$ and put
$\psi_6=\psi_3\eta_5$, so that
\[
\widehat{\Cl}=\{1,\eta_5\},\qquad \cR=\{\psi_3^{\pm1},\psi_6^{\pm1}\},\qquad
\cP=\{\psi_3,\psi_6\}.
\]
All four characters in $\cR$ have exact conductor $(2)$. Let $\rho$ be the
two-dimensional irreducible representation of $\Gal(CK/\Q)\cong S_3$, so that
$\zeta_C(s)=\zeta(s)L(s,\rho)$, and let $\chi$ be the quadratic character of $E/C$.

\begin{lemma}\label{lem:Artin35}
$\Ind_K^\Q\psi_3\cong\rho$ and $\Ind_K^\Q\psi_6\cong\rho\otimes\chi_5$. Consequently,
\begin{equation}\label{eq:Lh35}
L(h_{\psi_3},s)=L(s,\rho)=\frac{\zeta_C(s)}{\zeta(s)},\qquad
L(h_{\psi_6},s)=L(s,\rho\otimes\chi_5)=\frac{L(s,\chi)}{L(s,\chi_5)} .
\end{equation}
Both representations have conductor $140$ and root number $+1$, and
\begin{equation}\label{eq:L35values}
\begin{gathered}
L'(-1,\rho)=-12\,\zeta_C'(-1),\qquad L'(-1,\rho\otimes\chi_5)=-\tfrac52\,L'_{E/C}(-1,\chi),\\
L'(-1,V)=-\frac{140^{3/2}}{8\pi^3}L(2,V)
\end{gathered}
\end{equation}
for $V\in\{\rho,\rho\otimes\chi_5\}$.
\end{lemma}

\begin{proof}
Apply the standard Artin formalism \cite[Sections~2.2--2.3]{CogdellArtin}.
Since $\psi_3$ is a ring class character of order three, $\psi_3\circ c=\psi_3^{-1}$, so
$\Ind_K^\Q\psi_3$ is irreducible and factors through the $S_3$-quotient cut out by
$\ker\psi_3$; hence $\Ind_K^\Q\psi_3\cong\rho$. Because $\eta_5=\chi_5\circ\Nm_{K/\Q}$,
it follows that $\Ind_K^\Q(\psi_3\eta_5)\cong\rho\otimes\chi_5$. Moreover, $\sqrt5\notin C$, so
restricting $\chi_5$ to $\Gal(\overline\Q/C)$ gives $\chi$, and
$\Ind_C^\Q\chi=\chi_5\otimes\Ind_C^\Q1=\chi_5\oplus(\rho\otimes\chi_5)$. Since $\psi_3$
and $\psi_6$ are primitive of conductor $(2)$, $L(h_\psi,s)=L_K(\psi,s)=L(s,\Ind\psi)$
with conductor $|D_K|\Nm(2)=140$.

The root numbers of Dedekind zeta functions and of $\chi_5$ all equal $1$, so
$W(\rho)=W(\zeta_C)=1$ and $W(\rho\otimes\chi_5)=W(\zeta_E)W(\zeta_C)^{-1}W(\chi_5)^{-1}=1$.
Complex conjugation acts with eigenvalues $1,-1$ on both representations, giving gamma
factor $\Gamma_\C(s)=2(2\pi)^{-s}\Gamma(s)$; evaluating the functional equation at $s=-1$
gives the third relation in \eqref{eq:L35values}. The first two follow by differentiating
\eqref{eq:Lh35} at $s=-1$, using $\zeta(-1)=-\frac1{12}$,
\[
L(-1,\chi_5)=-\tfrac12B_{2,\chi_5}=-\tfrac12\cdot\tfrac{1-4-9+16}{5}=-\tfrac25,
\]
and the vanishing of $\zeta_C$, $L(s,\rho)$, $L(s,\chi)$, and $L(s,\rho\otimes\chi_5)$
at $s=-1$.
\end{proof}

\subsection{The evaluation}
Let $M_+=L'(g_\psi,0)$ and $M_-=L'(g_{\psi\eta_5},0)$ be the derivatives at $0$ of the
two weight-three CM newforms of level $35$ and nebentypus $\chi_{-35}$, and let
$\cU_{140}=\sum_{\chi'\in\cR}L'(g_{\psi\chi'},0)$ be the sum over the four newforms of
level $140$. In analogy with $\cD_d=-3H_{4d}$ when $h=1$, define
\begin{equation}\label{eq:D35}
\cD_{35}=-\frac3h\sum_{\rho'\in\cP}L'(h_{\rho'},-1)
=18\,\zeta_C'(-1)+\frac{15}4\,L'_{E/C}(-1,\chi),
\end{equation}
where the second equality follows from \eqref{eq:Lh35} and \eqref{eq:L35values}. As
before, $\Delta_a=L'(\chi_{-a},-1)$ denotes the separate Dirichlet $L$-derivative.

\begin{theorem}\label{thm:T35}
Theorem~D holds:
\[
n_4(s_{35})=\frac47(M_++M_-)+\frac4{21}\cU_{140}+\frac1{21}\Delta_{35}+\frac8{35}\Delta_7
+\frac{32}{7}\zeta_C'(-1)+\frac{20}{21}L'_{E/C}(-1,\chi).
\]
Equivalently, the last two terms equal $\frac{16}{63}\cD_{35}$.
\end{theorem}

\begin{proof}
Specialize Theorem~\ref{thm:GT} to $d=35$ and $h=2$. The first two coefficients are
$40/(dh)=4/7$ and $40/(3dh)=4/21$. For the third term, substituting \eqref{eq:genus35},
$\Delta_a=a^{3/2}L(\chi_{-a},2)/(4\pi)$, $\zeta(2)=\pi^2/6$, and
$L(\chi_5,2)=4\pi^2/(25\sqrt5)$ gives
\[
\frac{5\sqrt{35}}{2\pi^3}\,\zeta(2)L(\chi_{-35},2)=\frac1{21}\Delta_{35},\qquad
\frac{5\sqrt{35}}{2\pi^3}\,L(\chi_5,2)L(\chi_{-7},2)=\frac8{35}\Delta_7 .
\]
Finally, Lemma~\ref{lem:Artin35} turns the last term of Theorem~\ref{thm:GT} into
\[
-\frac{80}{3dh}\bigl(L'(-1,\rho)+L'(-1,\rho\otimes\chi_5)\bigr)
=\frac{8}{21}\Bigl(12\zeta_C'(-1)+\frac52L'_{E/C}(-1,\chi)\Bigr). \qedhere
\]
\end{proof}

\begin{remark}\label{rem:35general}
For general $d$, the final term of Theorem~A is $\frac{80}{9d}\cD_d$ with
$\cD_d=-\frac3h\sum_{\rho\in\cP}L'(h_\rho,-1)$, reducing to a single cubic zeta
derivative when $h=1$. At $d=35$, it splits into the cubic and sextic terms of
\eqref{eq:D35}, and the sextic piece encounters the Artin-factorization obstruction of
Proposition~\ref{prop:noKlein}.
\end{remark}

\begin{remark}
Numerically (see Table~\ref{tab:values} and Section~\ref{ssec:numerics}),
\begin{align*}
\zeta_C'(-1)&=0.51698\,09673\,79363\,15655\,74443\,8409\ldots,\\
L'_{E/C}(-1,\chi)&=2.79323\,00521\,11370\,37199\,44604\,5\ldots .
\end{align*}
\end{remark}

\section{Dilogarithms for \texorpdfstring{$d=35$}{d=35}}\label{sec:dilog}

The weight-two case of Zagier's conjecture is known \cite{Borel77,Zagier86,Suslin91}.
The issue here is the identification of particular short Bloch classes and their
exact rational coefficients. Table~\ref{tab:status} records the status of these
explicit identities; Sato's theorem supplies the ray-class component used below.

\begin{table}[ht]
\centering
\small
\renewcommand{\arraystretch}{1.22}
\begin{tabular}{@{}>{\raggedright\arraybackslash}p{0.37\textwidth}>{\raggedright\arraybackslash}p{0.19\textwidth}>{\raggedright\arraybackslash}p{0.36\textwidth}@{}}
\toprule
Statement & Status & Ingredients \\
\midrule
$\zeta_C'(-1)=\frac{D(\Xi_C)}{6\pi}$ (Theorem~\ref{thm:R1})
 & proved
 & Borel's volume formula, Selberg's lemma, the Dupont--Zickert volume cocycle;
   certificates K1--K4 \\
$A_3=-2\res\Xi_C$ (Proposition~\ref{prop:A3})
 & proved, using \cite{Sato26}
 & Sato's regulator identity, Galois descent, Borel--Suslin \\
$L'_{E/C}(-1,\chi)=\frac{6}{5\pi}D(\Xi_E)$, equivalently $A_6=-3\,\Xi_E$
(Conjecture~\ref{conj:R2})
 & conjecture
 & high-precision numerical agreement; exact reduction to one explicit element of
   $B(N)_\Q$ \\
Artin factorization through fields with at most one complex place
 & excluded (Proposition~\ref{prop:noKlein})
 & representation theory of $D_6$\\
\bottomrule
\end{tabular}
\caption{Status of the dilogarithmic statements for $d=35$.}
\label{tab:status}
\end{table}

\subsection{Notation}\label{ssec:dilog-notation}
Keep $C=\Q(a)$, $E=C(r)$ with $r^2=5$, and $N=E(t)$ with $t^2=-7$, the ray class field
of Lemma~\ref{lem:ray35}. Let $\sigma$ be the nontrivial automorphism of $E/C$, fixing $a$
and sending $r$ to $-r$. Fix the embedding of $C$ with $\Ima a>0$ and of $E$ with
$r=\sqrt5>0$. Because $E$ has real places, $K\not\subset E$, so there is a unique
embedding of $N$ extending this embedding of $E$ that maps $\sqrt{-35}$ to $i\sqrt{35}$.
View $N$ as a subfield of $\C$ via this embedding, and evaluate $D$ accordingly.

The subfield structure of $N$, together with the generators of
$\Gal(N/C)=\langle\kappa,\sigma\rangle$ used in Section~\ref{ssec:sato}, is shown below:
\[
\begin{array}{c@{\qquad}c@{\qquad}c}
 & N = C(\sqrt5,\sqrt{-7}) & \\[4pt]
\swarrow\ \scriptstyle\langle\kappa\rangle & \downarrow\ \scriptstyle\langle\kappa,\sigma\rangle & \searrow\ \scriptstyle C_6 \\[4pt]
E = C(\sqrt5) & C = \Q(a) & K = \Q(\sqrt{-35}) \\[4pt]
\downarrow\ \scriptstyle\langle\sigma\rangle & \downarrow\ \scriptstyle\deg 3 & \downarrow\ \scriptstyle\deg 2 \\[4pt]
C & \Q & \Q
\end{array}
\]

Recall the Bloch--Wigner dilogarithm $D(z)=\Ima\mathrm{Li}_2(z)+\arg(1-z)\log|z|$. For a
field $F$, $\mathcal P(F)$ denotes the rational pre-Bloch group,
$\delta[z]=z\wedge(1-z)\in\Lambda^2(F^\times)\otimes\Q$, and $B(F)_\Q=\ker\delta$. The convention for the five-term relation is
\[
R(x,y)=[x]-[y]+\Bigl[\frac yx\Bigr]-\Bigl[\frac{1-x^{-1}}{1-y^{-1}}\Bigr]+\Bigl[\frac{1-x}{1-y}\Bigr]=0.
\]
The standard rational symmetries $[z]=-[1-z]=-[1/z]$ follow after killing torsion
\cite[Sections~3.1--3.2]{Hutchinson13}. The cross-ratio convention is
\[
\crr(z_0,z_1,z_2,z_3)=\frac{(z_0-z_2)(z_1-z_3)}{(z_0-z_3)(z_1-z_2)}.
\]

The short elements $\Xi_C\in B(C)_\Q$ and $\Xi_E\in B(E)_\Q$ are those defined in
Section~\ref{sec:intro}. Verifying that they lie in the respective Bloch groups is a
finite check in $\Lambda^2$ of the $\{2,5\}$-units; for $\Xi_C$, this follows directly
from $1+a^5=(1-a)(1+a)^3$ (Section~\ref{ssec:exact}).

\subsection{The cubic identity}\label{ssec:R1}

\begin{theorem}\label{thm:R1}
$D(\Xi_C)=6\pi\,\zeta_C'(-1)$.
\end{theorem}

The proof follows Zagier's geometric method \cite{Zagier86}. Steps~1--2 identify
$v_0=\pi\zeta_C'(-1)$ as an arithmetic covolume. Steps~3--6 evaluate an explicit
$3$-cycle using integrality and a certified interval; Step~7 transfers the result
to $\Xi_C$ by five-term relations. Certificates K1--K4 supply the finite data.

The integrality step is indispensable, since an interval computation alone would give
only a numerical approximation. An analogous integrality bound is what is still missing
for Conjecture~\ref{conj:R2}, where Proposition~\ref{prop:noKlein} rules out repeating
this one-complex-place Artin-factorization argument.

\begin{proof}
Let $B=(-1,-1)_C$ and let $\mathcal O_B\subset B$ be the $\Z[a]$-order with the explicit
basis given in certificate K1; here $\cO_C=\Z[a]$.

\smallskip
\Needspace{90pt}
\noindent\textbf{Step 1 (maximal order).}\par
The discriminant of the reduced trace form of $\mathcal O_B$ is $\pm140^4\cdot2^2$. For
any order, $|\operatorname{disc}\mathcal O_B|=|d_C|^4\,\Nm(\Disc B)^2[\mathcal O_{\max}:\mathcal O_B]^2$
\cite{Voight21}.
Since $B$ is ramified at the real place of $C$, it must also be ramified at an odd
number of finite primes, so $\Nm(\Disc B)\ge2$. Hence $\mathcal O_B$ is maximal and
$\Disc B=\mathfrak p_2$, the unique prime above $2$ ($x^3+2x-2$ is Eisenstein at $2$).
In particular, $B$ is a division algebra.

\smallskip
\Needspace{90pt}
\noindent\textbf{Step 2 (covolume).}\par
Via a splitting $B\otimes_CL\cong M_2(L)$ over $L=C(i)$, the group
$\Gamma=\mathcal O_B^1/\{\pm1\}\subset\mathrm{PSL}_2(\C)$ is a cocompact Kleinian group.
By Borel's formula \cite{Borel81}, in the form \cite[Theorem~11.1.3]{MR03},
\[
v_0:=\covol(\Gamma)=\frac{|d_C|^{3/2}\zeta_C(2)}{(4\pi^2)^2}\prod_{\mathfrak p\mid\Disc B}(\Nm\mathfrak p-1)
=\frac{140^{3/2}\zeta_C(2)}{16\pi^4}.
\]
The functional equation for a field of signature $(1,1)$ gives
$\zeta_C'(-1)=|d_C|^{3/2}\zeta_C(2)/(16\pi^5)$, so $v_0=\pi\,\zeta_C'(-1)$.

\smallskip
\Needspace{90pt}
\noindent\textbf{Step 3 (an explicit $3$-cycle).}\par
Certificate K2 provides an explicit chain
$c=\sum_kw_k(1,g_{1,k},g_{2,k},g_{3,k})$ of $5520$ terms, supported on $41$ elements
of $\Gamma$ with weights in $\frac1{288}\Z$, whose alternation is a cycle in the
homogeneous bar complex of $\Gamma$ modulo $\Gamma$. This is verified in exact quaternion
arithmetic.

\smallskip
\Needspace{90pt}
\noindent\textbf{Step 4 (integrality).}\par
Let $\mathcal A(c)$ be the alternation of $c$, so that $6912\,\mathcal A(c)$ is integral
($6912=24\cdot288$). By Selberg's lemma \cite{Alperin87}, $\Gamma$ admits a
torsion-free subgroup $\Gamma'$ of finite index $m$. For the transfer
$\operatorname{tr}\colon H_3(\Gamma;\Z)\to H_3(\Gamma';\Z)$ and the inclusion $i$, $i_*\circ\operatorname{tr}=m$ \cite[III.9]{Brown82}. Since $M'=\HH^3/\Gamma'$ is
closed, $H_3(\Gamma';\Z)=\Z[M']$ and $\langle\mathrm{vol},[M']\rangle=m\,v_0$. Thus
$\langle\mathrm{vol},u\rangle\in\Z\,v_0$ for every $u\in H_3(\Gamma;\Z)$, which forces
$\langle\mathrm{vol},c\rangle\in\frac1{6912}\Z\,v_0$.

\smallskip
\Needspace{90pt}
\noindent\textbf{Step 5 (the volume cocycle).}\par
Define $V(z_0,\dots,z_3)=D(\crr(z_0,\dots,z_3))$ for distinct points in $\mathbb P^1(\C)$
and $V=0$ otherwise. Then $V$ is $\mathrm{PGL}_2$-invariant, alternating, and satisfies
the cocycle identity on all $5$-tuples, including degenerate ones. A prism homotopy
shows that the cocycles with base points $\infty$ and $y$ are cohomologous for every
$y\in\mathbb P^1(\C)$, and for generic $y$ this cocycle represents the volume class
\cite[Proposition~3.1]{Dupont87}, \cite[Section~1.3, Theorem~4.5]{DZ06}. Hence
$D(\beta)=\pm\langle\mathrm{vol},c\rangle$, where
\[
\beta=\sum_kw_k\bigl[\crr(\infty,g_{1,k}\infty,g_{2,k}\infty,g_{3,k}\infty)\bigr]\in\mathcal P(L)
\]
simplifies to $189$ distinct arguments after cancellation (the stabilizer of $\infty$ in
$\Gamma$ plays no role).

\smallskip
\Needspace{90pt}
\noindent\textbf{Step 6 (certified value).}\par
Ball arithmetic (certificate K3) encloses $\zeta_C(2)$ using its Euler product over
$p\le2\cdot10^6$ together with an explicit tail bound, yielding
$6912\,D(\beta)/v_0\in[-6912\pm0.011]$. Since Step~4 forces this quotient to be an
integer, it follows that $D(\beta)=-v_0$.

\smallskip
\Needspace{90pt}
\noindent\textbf{Step 7 (five-term certificate).}\par
In $\Q[L\setminus\{0,1\}]$, modulo the inversion relations $[z]+[1-z]$ and $[z]+[1/z]$,
certificate K4 establishes
\[
\Xi_C+6\beta=\sum_{k=1}^{450}c_k\,R(x_k,y_k),\qquad c_k\in\tfrac1{24}\Z,
\]
with explicit $x_k,y_k\in L$. Because $D$ vanishes on $R(x,y)$ and is odd under
$z\mapsto1-z$ and $z\mapsto1/z$, evaluating at the complex place of $L$ lying above the
chosen embedding of $C$ gives $D(\Xi_C)=-6D(\beta)=6v_0=6\pi\,\zeta_C'(-1)$.
\end{proof}

\begin{remark}\label{rem:R1-robust}
The Dupont--Sah cross-ratio convention \cite{DS82} used above is reciprocal to that of \cite{DZ06},
which accounts for the sign $D(\beta)=-v_0$. Note also that Steps~4 and~6 do not depend
on the exact scaling in Step~2: even if $\covol(\Gamma)$ were $2v_0$ or $v_0/2$, $-v_0$
would remain the unique lattice point inside the certified interval.

By \eqref{eq:L35values}, Theorem~\ref{thm:R1} is equivalent to
\[
L'(-1,\rho)=-\frac2\pi D(\Xi_C),\qquad
L(2,\rho)=\frac{16\pi^2}{140^{3/2}}D(\Xi_C).
\]
\end{remark}

\subsection{Sato's classes}\label{ssec:sato}
The enhanced conjecture for ideal classes originates in \cite{ZG00}; its ray-class
formulation is due to Sato \cite{Sato19}.
For modulus $\mathfrak f=(2)$, Sato's recent proof of the enhanced polylogarithm
conjecture in weight two provides classes $\xi_\nu\in B(N)_\Q$ indexed by
$\nu\in\Gal(N/K)$ \cite[Theorem~1.1]{Sato26}. Taking imaginary parts and converting his
convention $L_2=-D$, the regulator formula following \cite[Theorem~1.1]{Sato26} reads
\begin{equation}\label{eq:sato-reg}
D(\gamma^{-1}\xi_\nu)=\pi\,\zeta'_{(2)}(-1,\nu\gamma)\qquad(\nu,\gamma\in\Gal(N/K)),
\end{equation}
where $\gamma^{-1}\xi_\nu=\xi_{\nu\gamma}$ and $\zeta_{(2)}(s,\nu)$ is the partial zeta
function of the ray class corresponding to $\nu$. Sato's theorem is used as a cited
result without reproving it here.

Let $c$ denote complex conjugation; since $a\notin\R$, $c$ does not fix $C$. The
subgroup $\Gal(N/C)$ has order four and contains the unique involution $g^3$ of the
cyclic group $\Gal(N/K)=\langle g\rangle$ together with two elements $g^kc$ and
$g^{k+3}c$. Replacing $g$ by $g^{-1}$ if necessary, arrange that $\kappa:=g^2c$
fixes $C$. Put $\sigma:=g^3$ and $\xi_j=\xi_{g^j}$ ($j\bmod6$), so that
$\Gal(N/C)=\langle \kappa,\sigma\rangle$. Since $g$ acts nontrivially on the genus field
and $c$ fixes $\sqrt5$, the element $\kappa$ fixes $\sqrt5$; hence $E=N^{\langle\kappa\rangle}$,
and $\sigma|_E$ coincides with the automorphism $\sigma$ defined earlier.

\begin{lemma}\label{lem:xi-action}
$g\xi_j=\xi_{j-1}$ and $c\,\xi_j=-\xi_{-j}$. Consequently,
$\kappa\xi_j=-\xi_{-j-2}$ and $\sigma\xi_j=\xi_{j-3}$.
\end{lemma}

\begin{proof}
The first relation is $\gamma^{-1}\xi_\nu=\xi_{\nu\gamma}$ with $\gamma=g^{-1}$. For the
second, the regulator map $x\mapsto(D(\iota x))_{\iota}$ over all embeddings $\iota$ of
$N$ is injective on $B(N)_\Q$ by the theorems of Borel and Suslin
\cite{Borel77,Suslin91}. These embeddings are $\gamma$ and $c\gamma$ with
$\gamma\in\Gal(N/K)$, composed with the fixed embedding. Since $\Gal(N/\Q)$ is
generalized dihedral, $c\gamma c=\gamma^{-1}$ and $D(\bar z)=-D(z)$. Moreover,
$\zeta_{(2)}(s,\nu)=\zeta_{(2)}(s,\nu^{-1})$ because the conjugate of an ideal lies in
the inverse ray class and has the same norm. Applying \eqref{eq:sato-reg} yields
\[
D(\gamma\,c\,\xi_\nu)=-D(\gamma^{-1}\xi_\nu)=-\pi\zeta'_{(2)}(-1,\nu\gamma)
=D(\gamma(-\xi_{\nu^{-1}})),
\]
and similarly at the embeddings $c\gamma$, whence $c\xi_\nu=-\xi_{\nu^{-1}}$.
\end{proof}

Set
\[
A_3=\xi_0-\xi_4+\xi_3-\xi_1,\qquad A_6=\xi_0-\xi_4-\xi_3+\xi_1,\qquad
\Gamma_{35}=\xi_0-\xi_4=\tfrac12(A_3+A_6).
\]

\begin{proposition}\label{prop:satoL}
$D(A_3)=\pi L'(-1,\rho)=-12\pi\,\zeta_C'(-1)$ and
$D(A_6)=\pi L'(-1,\rho\otimes\chi_5)=-\frac{5\pi}2L'_{E/C}(-1,\chi)$.
In particular, $D(\Gamma_{35})=-\frac\pi3\cD_{35}$.
\end{proposition}

\begin{proof}
Let $\chi_k$ be the character of $\Gal(N/K)$ defined by $\chi_k(g)=e^{\pi ik/3}$. By
character orthogonality, $\zeta_{(2)}(s,g^j)=\frac16\sum_{k=0}^5\overline{\chi_k(g^j)}L_{(2)}(s,\chi_k)$,
so for any combination $A=\sum_jc_j\xi_j$, \eqref{eq:sato-reg} with $\gamma=1$ gives
\[
D(A)=\frac\pi6\sum_{k\bmod 6}\Bigl(\sum_jc_je^{-\pi ijk/3}\Bigr)L'_{(2)}(-1,\chi_k).
\]
For $A_3$, the inner character sums are $0,0,3+i\sqrt3,0,3-i\sqrt3,0$ for $k=0,\dots,5$;
for $A_6$, they are $0,3-i\sqrt3,0,0,0,3+i\sqrt3$. The characters $\chi_{\pm2}$ have
order three and $\chi_{\pm1}$ have order six; all four are primitive of conductor $(2)$,
and inverse characters share the same $L$-function. By Lemma~\ref{lem:Artin35}, these
$L$-functions are $L(s,\rho)$ and $L(s,\rho\otimes\chi_5)$, which establishes the first
two formulas, while the third follows from \eqref{eq:D35}.
\end{proof}

This identity is independent of the choice of generator $g$. Using Sato's theorem \cite{Sato26}, Proposition~\ref{prop:satoL} and Theorem~D give
\begin{equation}\label{eq:n4-Gamma}
n_4(s_{35})=\frac47(M_++M_-)+\frac4{21}\cU_{140}+\frac1{21}\Delta_{35}+\frac8{35}\Delta_7
-\frac{16}{21\pi}D(\Gamma_{35}).
\end{equation}
However, $\Gamma_{35}$ is a long element whose explicit representatives involve
thousands of symbols in auxiliary extensions of $N$.

\begin{proposition}\label{prop:A3}
$A_3=-2\res_{N/C}\Xi_C$ in $B(N)_\Q$.
\end{proposition}

\begin{proof}
By Lemma~\ref{lem:xi-action}, $\kappa A_3=\sigma A_3=A_3$, so rational Galois descent
for Bloch groups \cite[Section~3]{BRV02} gives $A_3=\res_{N/C}\alpha$ for some
$\alpha\in B(C)_\Q$. The space $B(C)_\Q$ has dimension $r_2(C)=1$, and its
regulator is injective \cite{Borel77,Suslin91}. Since $D(\Xi_C)\ne0$, it follows
that $\alpha=\lambda\Xi_C$ with
$\lambda=D(A_3)/D(\Xi_C)=-12\pi\zeta_C'(-1)/(6\pi\zeta_C'(-1))=-2$ by
Proposition~\ref{prop:satoL} and Theorem~\ref{thm:R1}.
\end{proof}

\Needspace{155pt}
\begin{corollary}\label{cor:A6}
The following statements are equivalent:
\begin{enumerate}
\item $L'_{E/C}(-1,\chi)=\frac6{5\pi}D(\Xi_E)$;
\item $A_6=-3\,\Xi_E$ in $B(E)_\Q$;
\item $\Gamma_{35}=-\res_{N/C}\Xi_C-\tfrac32\,\Xi_E$ in $B(N)_\Q$.
\end{enumerate}
\end{corollary}

\begin{proof}
By Lemma~\ref{lem:xi-action}, $\kappa A_6=A_6$ and $\sigma A_6=-A_6$. Thus both $A_6$
and $\Xi_E$ lie in $B(E)_\Q^{\sigma=-1}$, which has dimension $r_2(E)-r_2(C)=1$. Since
$D(\Xi_E)\ne0$ by the interval bound in Section~\ref{ssec:enclosures}, write $A_6=\lambda\Xi_E$ with $\lambda=D(A_6)/D(\Xi_E)$, and
Proposition~\ref{prop:satoL} shows that $\lambda=-3$ if and only if~(1) holds, giving
(1)$\Leftrightarrow$(2). Equivalence with~(3) follows from
$\Gamma_{35}=\tfrac12(A_3+A_6)$ and Proposition~\ref{prop:A3}.
\end{proof}

\subsection{The sextic identity}\label{ssec:R2}

\begin{conjecture}\label{conj:R2}
$L'_{E/C}(-1,\chi)=\frac{6}{5\pi}D(\Xi_E)$. Equivalently,
\[
L'(-1,\rho\otimes\chi_5)=-\frac3\pi D(\Xi_E),\qquad
L(2,\rho\otimes\chi_5)=\frac{24\pi^2}{140^{3/2}}D(\Xi_E),\qquad
\cD_{35}=\frac3\pi D(\Xi_C)+\frac9{2\pi}D(\Xi_E).
\]
\end{conjecture}

The equivalences follow immediately from \eqref{eq:L35values}, \eqref{eq:D35}, and
Theorem~\ref{thm:R1}. If Conjecture~\ref{conj:R2} holds, then
\[
n_4(s_{35})=\frac47(M_++M_-)+\frac4{21}\cU_{140}+\frac1{21}\Delta_{35}+\frac8{35}\Delta_7
+\frac{16}{21\pi}D(\Xi_C)+\frac{8}{7\pi}D(\Xi_E).
\]

Because $E$ has two complex places, the hyperbolic volume argument of
Section~\ref{ssec:R1} does not apply directly. The next proposition rules out a natural
workaround: realizing $\rho\otimes\chi_5$ inside the permutation representation of any
number field with at most one complex place.

\begin{proposition}\label{prop:noKlein}
Let $F$ be a number field with at most one complex place. Then $\rho\otimes\chi_5$ does
not occur in the permutation representation $\Ind_F^\Q1$. Consequently,
$L(s,\rho\otimes\chi_5)$ is not a product of rational powers of Dedekind zeta functions of
such fields.
\end{proposition}

\begin{proof}
Let $M$ be a finite Galois extension containing $N$ and $F$, and put $G=\Gal(M/\Q)$,
$H=\Gal(M/F)$, and $X=G/H$, so that $\Ind_F^\Q1=\C[X]$. Complex conjugation $c\in G$
acts on $X$ just as it acts on the complex embeddings of $F$, fixing the real
embeddings and pairing the complex ones into $r_2(F)$ transpositions. Hence the
$(-1)$-eigenspace of $c$ on $\C[X]$ has dimension $r_2(F)$.

\smallskip
\Needspace{90pt}
\noindent\textbf{Case $r_2(F)=0$.}\par
On $\rho\otimes\chi_5$, complex conjugation has eigenvalues $1,-1$, since it acts on the
roots of $x^3+2x-2$ as a transposition while $\chi_5(c)=1$. Thus $r_2(F)=0$ immediately
excludes $\rho\otimes\chi_5\subset\C[X]$.

\smallskip
\Needspace{90pt}
\noindent\textbf{Case $r_2(F)=1$.}\par
Suppose $r_2(F)=1$ and $\rho\otimes\chi_5\subset\C[X]$. Then the multiplicity is one, and
$c$ acts trivially on the $G$-stable complement $V'$ of this copy. The kernel of $G$ on
$V'$ is normal and contains $c$, so it contains the normal closure $J$ of $c$ in $G$.
The image of $J$ in $\Gal(N/\Q)$ is $\Gal(N/\Q(\sqrt5))\cong S_3$, because $\Q(\sqrt5)$
is the maximal totally real subfield of $N$, and on this subgroup $\rho\otimes\chi_5$
restricts to the irreducible representation $\rho$.

Thus the permutation character of $J$ on $X$ is $(|X|-2)\cdot1+\rho|_J$, which has
$|X|-2$ orbits whose nontrivial parts sum to the irreducible two-dimensional character
$\rho|_J$. Consequently, $J$ must have a single orbit $O$ of size three and fix all
remaining points of $X$. Since $G$ normalizes $J$, it preserves both $O$ and the
fixed-point set, so transitivity forces $X=O$ and $[F:\Q]=3$. However, the permutation
representation of a cubic field factors through a group of order at most $6$, whereas
$\rho\otimes\chi_5$ is faithful on $\Gal(N/\Q)\cong S_3\times C_2$ of order $12$, a
contradiction.

Finally, since Artin $L$-functions of distinct irreducible characters are
multiplicatively independent and $\rho\otimes\chi_5$ has multiplicity zero in every
$\Ind_F^\Q1$ under consideration, no rational-power factorization into their Dedekind
zeta functions can exist.
\end{proof}

By Borel's covolume formula \cite{Borel81}, the volume of an arithmetic hyperbolic
$3$-orbifold is a rational multiple of $|d_k|^{3/2}\zeta_k(2)\pi^{-2[k:\Q]+2}$ where the
invariant trace field $k$ has a single complex place. Proposition~\ref{prop:noKlein}
therefore rules out deriving Conjecture~\ref{conj:R2} from an Artin factorization of
$L(s,\rho\otimes\chi_5)$ into Dedekind zeta functions of such fields, though it does
not exclude identities at the single argument $s=2$ or period constructions over $E$.

\begin{remark}[Evidence]\label{rem:R2-evidence}
At $130$-digit working precision,
\[
\bigl|D(\Xi_E)-(5\pi/6)L'_{E/C}(-1,\chi)\bigr|\approx4.51\cdot10^{-130},
\]
with $L(2,\rho\otimes\chi_5)$ computed from Epstein zeta functions of discriminant
$-140$. This floating-point comparison is not a rigorous interval enclosure and
does not prove the Bloch-class identity in Corollary~\ref{cor:A6}.
\end{remark}

An explicit five-term certificate for $A_6+3\Xi_E=0$, a direct evaluation, or a
denominator bound combined with a sufficiently tight certified interval would close
Conjecture~\ref{conj:R2}. Rational proportionality alone does not determine the
coefficient; compare \cite{BRV05}.

\section{Computations and certificates}\label{sec:computations}

Finite computations used in the proofs consist of exact identities and certified
interval bounds, supplementing the analytic and geometric arguments.
Floating-point comparisons provide numerical evidence only.

\subsection{Exact algebraic checks}\label{ssec:exact}
\begin{itemize}
\item The parameter resultants, irreducibility witnesses, and root brackets of
Section~\ref{sec:parameters}, as well as the field identifications, discriminants, and
Dedekind maximality witnesses of Section~\ref{sec:classone}, were verified independently
in exact arithmetic over $\Z$.
\item In $C$, the identity $1+a^5=(1-a)(1+a)^3$ immediately gives
$\delta\Xi_C=5\,a\wedge(1+a^5)-5\,a\wedge(1-a)-15\,a\wedge(1+a)=0$ modulo torsion.
\item Exact decompositions into $\{2,5\}$-units of $E$ give
$\delta\Xi_E=0$ in $\Lambda^2$ and also confirm $\delta\Xi_C=0$.
\end{itemize}

\subsection{Hecke coefficients}\label{ssec:hecke}
Weight-three newforms enter the formulas through $L'(g,0)$ (equivalently $L(g,3)$), and
weight-one forms through $L'(h,-1)$ (equivalently $L(h,2)$). Table~\ref{tab:levels}
collects the primitive levels, nebentypus characters, and scalar terms appearing in
Theorems~B--D.

\begin{table}[ht]
\centering
\small
\renewcommand{\arraystretch}{1.2}
\setlength{\tabcolsep}{9pt}
\begin{tabular}{@{}rcccl@{}}
\toprule
$d$ & Weight $3$ levels & Weight $1$ level & Nebentypus & Scalar values\\
\midrule
$11$ & $11,\,44$ & $44$ & $\chi_{-11}$ & $\Delta_{11},\ \zeta'_{F_{11}}(-1)$\\
$19$ & $19,\,76$ & $76$ & $\chi_{-19}$ & $\Delta_{19},\ \zeta'_{F_{19}}(-1)$\\
$27$ & $12,\,27,\,108$ & $108$ & $\chi_{-3}$ & $\Delta_3,\ \zeta'_{F_{27}}(-1)$\\
$43$ & $43,\,172$ & $172$ & $\chi_{-43}$ & $\Delta_{43},\ \zeta'_{F_{43}}(-1)$\\
$67$ & $67,\,268$ & $268$ & $\chi_{-67}$ & $\Delta_{67},\ \zeta'_{F_{67}}(-1)$\\
$163$ & $163,\,652$ & $652$ & $\chi_{-163}$ & $\Delta_{163},\ \zeta'_{F_{163}}(-1)$\\
$35$ & $35,\,140$ & $140$ & $\chi_{-35}$ & $\Delta_{35},\Delta_7,\zeta'_C(-1),L'_{E/C}(-1,\chi)$\\
\bottomrule
\end{tabular}
\caption{Primitive levels and special values. The nebentypus applies to both
weights, with the quadratic character viewed modulo the indicated level.
There are two weight-one forms at level $140$.}
\label{tab:levels}
\end{table}

Although Proposition~\ref{prop:decomp}, Lemma~\ref{lem:traces}, and
Proposition~\ref{prop:F-lat} hold for all Fourier coefficients, two independent
constructions of the $q$-expansions were also compared through $q^{32768}$ (one from local
Euler factors and primary prime generators, the other by direct enumeration of lattice
vectors). All $655\,360$ positive-index coefficients agree, and negative controls
detect omitted oldforms in \eqref{eq:Flat}: deleting $3g_{12}(3z)$ shifts the $q^3$
coefficient from $0$ to $-3$, while deleting $27g_{12}(9z)$ shifts the $q^9$
coefficient from $27$ to $0$.

\subsection{Certified interval bounds}\label{ssec:enclosures}
For the six cases of Theorems~B and~C, rigorous enclosures of $n_4(s_d)$ of
width below $10^{-90}$ were computed using outward-rounded dyadic interval arithmetic with denominator
$2^{360}$. Integer bisection isolates the root $s=-A$ of $Q_{j_d}$, and the exterior
series of Lemma~\ref{lem:local} is summed to $45$ terms with tail bound
$\varrho^{46}/(46(1-\varrho))$ ($\varrho=256/A<1$). All $14\,000$ endpoint checks
agree with exact rational arithmetic. Note that these intervals enclose the left-hand
sides $n_4(s_d)$ only, not the right-hand sides or the residual in Conjecture~\ref{conj:R2}.

A separate Arb evaluation at the embedding fixed in Section~\ref{ssec:dilog-notation}
gives $731/100<D(\Xi_E)<732/100$. This proves the nonvanishing used in
Corollary~\ref{cor:A6}, but establishes no identity with an $L$-value.

\subsection{Numerical comparison}\label{ssec:numerics}
At $60$ and $130$ decimal digits of working precision, a separate script compares:
(1) $s_d$ from the eta quotient \eqref{eq:t-s4} and $n_4(s_d)$ from the exterior
series; (2) $n_4(s_d)$ from the lattice formula \eqref{eq:lattice} via one-dimensional
Poisson summation; and (3) the right-hand sides of Theorems~B--D, where the
order-three constants are evaluated via Epstein zeta functions restricted to values
coprime to the conductor (independently of the coset decomposition in
Section~\ref{sec:general}). The weight-three values reuse lattice sums, so these
comparisons are only partly independent. The observed $130$-digit residuals for
Theorems~B--D are all below $2.2\cdot10^{-129}$; they are not interval enclosures.

\begin{table}[ht]
\centering
\small
\setlength{\tabcolsep}{9pt}
\renewcommand{\arraystretch}{1.18}
\begin{tabular}{@{}rr@{\,}lll@{}}
\toprule
$d$ & \multicolumn{2}{c}{$s_d$} & $n_4(s_d)$ & $\zeta'_{F}(-1)$\\
\midrule
$11$ & $-3.34022\,73463\,4$ & $\cdot10^{4}$ & $10.41709\,66323\,51$ & $0.08428\,03218\,23$ \\
$19$ & $-8.85375\,78261\,7$ & $\cdot10^{5}$ & $13.69379\,45525\,21$ & $0.21063\,04063\,82$ \\
$27$ & $-1.22886\,39984\,3$ & $\cdot10^{7}$ & $16.32418\,77680\,86$ & $0.36737\,41284\,28$ \\
$43$ & $-8.84736\,64000\,0$ & $\cdot10^{8}$ & $20.60080\,06038\,91$ & $0.75372\,92705\,36$ \\
$67$ & $-1.47197\,95264\,0$ & $\cdot10^{11}$ & $25.71504\,41346\,11$ & $1.48035\,55953\,10$ \\
$163$ & $-2.62537\,41264\,1$ & $\cdot10^{17}$ & $40.10916\,99911\,32$ & $5.65994\,47893\,13$ \\
$35$ & $-1.17966\,57776\,5$ & $\cdot10^{8}$ & $18.58591\,21064\,79$ & $0.51698\,09673\,79$ \\
\bottomrule
\end{tabular}
\caption{Values from the $60$-digit run, truncated for display: $s_d$ from the eta
quotient, $n_4(s_d)$ from the exterior series, and $\zeta'_F(-1)$ from the scalar
Epstein calculation ($F=F_d$ for $h=1$, $F=\Q(\sqrt[3]2)$ for $d=27$, and $F=C$ for
$d=35$). The three evaluations of $n_4(s_d)$ agree to within $2\cdot10^{-59}$;
rigorous interval enclosures are separate (Section~\ref{ssec:enclosures}).}
\label{tab:values}
\end{table}

At $d=35$, the same computation compares $D(\Xi_C)$ and $D(\Xi_E)$ against the
corresponding $L$-values; only the cubic identity is proved (Theorem~\ref{thm:R1}). As a
further normalization check at $d=27$, replacing $31M_{12}$ by $30M_{12}$ in Theorem~C
shifts the right-hand side by $40M_{12}/81\approx0.1489$.

\Needspace{170pt}
\subsection{Certificates for Theorem E}\label{ssec:certs}
Together with the written integrality argument in Section~\ref{ssec:R1}, the four
certificates K1--K4 form the proof-critical inputs for Theorem~\ref{thm:R1}:

\begin{description}
\item[K1] The basis of $\mathcal O_B\subset(-1,-1)_C$, its reduced trace form, and its
discriminant $\pm140^4\cdot2^2$, verified in exact arithmetic.
\item[K2] The chain $c$ of $5520$ terms on $41$ elements of $\Gamma$ with weights in
$\frac1{288}\Z$, whose cycle condition is verified in exact quaternion arithmetic.
\item[K3] The Arb ball-arithmetic enclosure \cite{Arb}
$6912\,D(\beta)/v_0\in[-6912\pm0.011]$, bounding $\zeta_C(2)$ by its Euler product
over $p\le2\cdot10^6$ with tail factor in $[1,e^{3/X}]$ ($X=2\cdot10^6$).
\item[K4] The $450$ five-term relations with coefficients in $\frac1{24}\Z$ expressing
$\Xi_C+6\beta$, verified modulo inversion relations by exact arithmetic in $L=C(i)$.
\end{description}
Independent readers rechecked all three exact rational certificates as well as the Arb
interval enclosure \cite{Arb}.

\subsection{Evidence for Conjecture F}
The ancillary files verify the closedness and seven-symbol form of $\Xi_E$ in exact
arithmetic, certify its nonzero regulator, and reproduce the floating-point
comparison in Remark~\ref{rem:R2-evidence}. They do not prove $A_6=-3\Xi_E$.

\subsection{Software and assistance}
The computations were carried out in Python (using SymPy, mpmath, and python-flint/Arb
\cite{Arb}) and PARI/GP \cite{PARI}. AI-based tools assisted with parts of the
computations and text preparation. The author takes responsibility for the contents
of this article.\par

The author pursues independent mathematical research alongside his main professional
work. Comments and corrections are welcome.

\microtypesetup{expansion=false}

\end{document}